\documentclass[12pt,a4paper,twoside]{amsart}
\usepackage{a4,amscd,amssymb,amsmath,amsthm,eucal,fancyheadings,xypic}
\usepackage[all]{xy}
\usepackage[english]{babel}
\numberwithin{equation}{section}
\newtheorem{thm}{Theorem}[section]

\newtheorem{lemma}[thm]{Lemma}
\newtheorem{prop}[thm]{Proposition}
\newtheorem{cor}[thm]{Corollary}

\newtheorem{ass}[]{Assumption}

\newcommand{\z}{\zeta}
\newcommand{\al}{\alpha}
\newcommand{\be}{\beta}
\newcommand{\g}{\gamma}
\newcommand{\ga}{\gamma}
\newcommand{\e}{\epsilon}
\newcommand{\la}{\lambda}
\newcommand{\Z}{\mathbb Z}
\newcommand{\Q}{\mathbb Q}

\newcommand{\F}{\mathbb F}
\newcommand{\V}{\mathcal V}
\newcommand{\m}{\mbox{ mod }}

\DeclareMathOperator{\pic}{Pic}
\DeclareMathOperator{\cl}{Cl}
\DeclareMathOperator{\kernel}{ker}
\DeclareMathOperator{\trace}{Tr}
\DeclareMathOperator{\Image}{Im}

\DeclareMathOperator{\cha}{Char}
\DeclareMathOperator{\Gal}{Gal}
\DeclareMathOperator{\aut}{Aut}

\def\proofname{\bf Proof.}
\begin{document}
\title[The Kervaire-Murthy Conjectures]{Fine Structure of Class Groups
$\cl^{(p)}\Q(\z_n)$
%of Prime Power Cyclic Groups
and the Kervaire-Murthy Conjectures}
\author{Ola Helenius}
\address{Department of Mathematics, Chalmers University of Technology
  and G\"{o}teborg University, SE-41296 G\"{o}teborg, Sweden}
%\curraddr{}
\email{olahe@math.chalmers.se, astolin@math.chalmers.se}

\author{Alexander Stolin}
%\address{Department of Mathematics, Chalmers University of Technology and G\"{o}teborg University, SE-41296 G\"{o}teborg, Sweden}
%\curraddr{}
%\email{astolin@math.chalmers.se}

%    General info
\date{}
%
%\dedicatory{}
%
\begin{abstract}
In 1977 Kervaire and Murthy presented three conjectures regarding $K_0 \Z C_{p^n}$, where $C_{p^n}$ is
the cyclic group of order $p^n$ and $p$ is a semi-regular prime that is $p$ does not divide
$h^+$ (regular $p$ does not divide the class number $h=h^+ h^-$).
The Mayer-Vietoris exact sequence provides the following short exact sequence
$$0\to V_n\to \pic \Z C_{p^n}\to \cl \Q (\z_{n-1})\times   \pic \Z C_{p^{n-1}}\to 0$$
Here $\z_{n-1}$ is a primitive $p^n$-th root of unity.
The group $V_n$ that injects into
$\pic \Z C_{p^n}\cong\tilde{K}_0\Z C_{p^n} $, is a canonical quotient of an in some sense simpler
group
$\V_n$. Both groups split in a ``positive'' and ``negative'' part. While $V_n^-$ is well understood there is
still no complete information on  $V_n^+$.
Kervaire and Murthy showed that
$K_0 \Z C_{p^n}$ and $V_n$ are tightly connected to class groups
of cyclotomic fields.
They also conjectured that $V_n^+ \cong(\Z/p^n\Z)^{r(p)}$,
where
$r(p)$ is the index of regularity of the prime $p$ and that $\V_n^+\cong V_n^+$, and moreover,
$\cha \V_n^+\cong \cl^{(p)} \Q (\z_{n-1})$, the $p$-part of the class group.

Under an extra assumption on
the prime $p$, Ullom proved in 1978
%in \cite{U2}
that
$V_n^+\cong (\Z/p^n\Z)^{r(p)}\oplus (\Z/p^{n-1}\Z)^{\la-r(p)}$, where $\la$ is one of the Iwasawa invariants.
Hence Kervaire and Murthys first conjecture holds only when
$\la=r(p)$.

In the present paper we calculate $\V_n^+$  and prove that $\cha \V_n^+\cong \cl^{(p)} \Q(\z_{n-1})$ for all
semi-regular primes which also gives us the structure of $\cl^{(p)} \Q(\z_{n-1})$ as an abelian group.
We also prove that under the same condition Ullom used, conjecture two always holds, that is
$\V_n^+\cong V_n^+$.
Under the assumption $\la=r(p)$ we construct a
special basis for a ring closely related to $\Z C_{p^n}$, consisting of units from a number field. This basis
is used to prove
that $\V_n^+\cong V_n^+$ in this case and it also follows that the Iwasawa invariant $\nu$ equals $r(p)$.
Moreover we conclude that $\la=r(p)$ is equivalent to that all three Kervaire and Murthy conjectures hold.

\end{abstract}
\maketitle

\section{Introduction}

Let $p$ be an odd prime,  $C_{p^n}$ denote the cyclic group of
order $p^n$ and let $\z_n$ be a primitive $p^{n+1}$th root of
unity. In this paper we work on the problem of finding $\pic \Z
C_{p^n}$. Our methods also lead to the calculation of the $p$-part
of the ideal class group of $\Z[\z_n]$. Calculating Picard groups
for a group ring like the one above is equivalent to calculating
$K_0$ groups. Finding $K_0 \Z G$ for various groups $G$ was
mentioned by R.G. Swan at his talk at the International Congress
of Mathematicians in Nice 1970 as one of the important problems in
algebraic K-theory. Of course, the reasons for this is are
applications in topology. However, calculating $K_0 (\Z G)$ seems
to be pretty hard and even to this date there are no general
results. Even when we restrict ourselves to $G=C_{p^n}$ no general
explicit formulas are known. Several people have worked on this,
though. Kervaire and Murthy presented in \cite{K-M} an approach
based on the pullback
\begin{equation}\label{eq:KM-pullback}
\xymatrix@=40pt{
\Z C_{p^{n+1}}  \ar[r] \ar[d]  & \Z[\z_n]  \ar[d]\\
\Z C_{p^{n}} \ar[r]         &  \F_p[x]/(x^{p^n}-1)=:R_n
}
\end{equation}
%which is a variant of \ref{eq:idealpullback}.
The $(*,\pic)$-Mayer-Vietoris exact sequence associated
to this pullback reads
\[
(\Z C_{p^{n}})^*\times \Z[\z_n]^* \stackrel{j}{\rightarrow}
R_n^*\to \pic \Z C_{p^{n+1}} \to \pic \Z C_{p^{n}}\times \pic \Z[\z_n] \to
\pic R_n
\]
Following Kervaire and Murthy, we observe that Picard groups of local rings
are trivial, that the Picard group of a Dedekind ring equals the class group of the same ring
and then define $V_n$ as the co-kernel of the map $j$ in the sequence above. Then we get
\begin{equation}\label{ZC-MV}
0\to V_n \to \pic \Z C_{p^{n+1}} \to \pic \Z C_{p^{n}}\times \cl \Q(\z_n) \to 0.
\end{equation}
Kervaire and Murthy set out to calculate $V_n$ and their approach is based on the fact that
all rings involved can be acted upon by the Galois group $G_n:=\Gal(\Q(\z_n)/\Q)$. If $s\in G_n$,
let $s(\z_n)=\z_n^{\kappa(s)}$. If we represent the rings in the pullback as residue class rings of
polynomials in the indeterminate $X$, the action is generated by $s(X)=X^{\kappa(s)}$ for all
involved rings. $G_n$ becomes a group of automorphisms of $\Z C_{p^{n+1}}$, $\Z C_{p^{n}}$
and $R_n$. The maps
in the pullback above commutes with the action of $G_n$ and the exact sequence becomes a sequence
of $G_n$-modules. In particular, complex conjugation, which we denote by $c$, belongs to
$G_n$ and  $c(X)=X^{-1}$. When $M$ is a multiplicative $G_n$-module, like the group of units of
one of the rings in the pullback, we let $M^+$ denote the subgroup of elements $v\in M$ such that
$c(v)=v$ and $M^-$ denote the subgroup of elements such that
$c(v)=v^{-1}$. $V_n$ is a finite abelian group of odd order and hence we have that
$V_n=V_n^+ \times V_n^-$. The main result in Kervaire and Murthy's article is the following theorem
\begin{thm}[Kervaire and Murthy]
\[
V_n^-\cong \prod_{\nu =1}^{n-1} (\Z/p^{\nu}\Z)^{\frac{(p-1)^2 p^{n-\nu-1}}{2}}
\]
and when $p$ is semi-regular, there exists a canonical injection
\[
\cha V_n^+ \to \cl^{(p)} \Q(\z_{n-1}),
\]
where $\cl^{(p)} \Q(\z_{n-1})$ is
the $p$-primary component of the ideal class group of $\Q(\z_{n-1})$.
\end{thm}
The calculation of $V_n^-$
is straightforward.
%Kervaire and Murthy first uses an extension to $\Z C_{p^n}$
%of Kummer's well known result that a unit in $\Z[\z_n]$ can be written as a product of a real unit and
%a power of $\z_n$ to show that $R_n^*/(X_n (R_n^*)^+$ and then the result follows from
%an analysis of the group $R_n^*$.
Finding the information on $ V_n^+$ turns out to be much harder. Kervaire and Murthy instead
proves the result above with $V_n^+$ replaced by
the $+$-part of
\[
\V_n:=\frac{R_n^*}{j(\Z[\z_n]^*)},
\]
that is, constructs a canonical injection
\begin{equation}\label{eq:KMinj}
\cha \V_n^+ \to \cl^{(p)} \Q(\z_{n-1})
\end{equation}
Then, since $V_n^+$ is a canonical quotient of $\V_n^+$,
\ref{eq:KMinj}
%this injection
extends to an injection
\[
\cha V_n^+ \to \cl^{(p)} \Q(\z_{n-1})
\]
via the canonical injection
\[
\cha V_n^+ \to \cha \V_n^+.
\]

The injection \ref{eq:KMinj}
is actually a composition of the Artin map in class field theory and a canonical injection
from Iwasawa theory. The actual proof is mainly based on class field theory.

Let $r(p)$ be the index of regularity of $p$, that is the number of Bernoulli numbers
$B_2, B_4,\hdots,B_{p-3}$ with numerators (in reduced form) divisible by $p$.
%Let $\cl^{(p)} \Z(\z_{n-1})$ denote the $p$-component of the ideal class group of
%$\Z(\z_{n-1})$.
Kervaire and Murthy formulate the following conjectures. For semi-regular primes:
\begin{eqnarray}
V_n^+ &=& \V_n^+ \label{eq:conj1}\\
\cha \V_n^+ &=& \cl^{(p)} \Q(\z_{n-1}) \label{eq:conj3}\\
\cha V_n^+ &\cong& \big(\frac{\Z}{p^n\Z}\big)^{r(p)},\label{eq:conj2}
\end{eqnarray}

When $p$ is a regular prime it is known that $\cl^{(p)} \Q(\z_{n-1})$ is trivial and hence
$V_n=V_n^-$ is determined completely in \cite{K-M}.

In \cite{U2}, Stephen Ullom uses Iwasawa theory and studies the action of $\aut C_{p^n}$ on
$\pic \Z C_{p^n}$.
He proves in that under a certain extra assumption on $p$,
%condition on the
%Iwasawa invariants associated to the semi-regular prime $p$,
%
%$V_n^+ \cong \big(\frac{\Z}{p^n\Z}\big)^r$
the first of Kervaire and Murthy's conjectures hold exactly when the Iwasawa invariant $\la$
associated to $p$ equals $r(p)$.
%conjecture~\ref{eq:conj2} holds.
%%%%%%%%%%%%
%More explicitly the assumption is the following.
%\begin{ass}\label{ass:Ullomcond}
%The Iwasawa
%invariants $\la_{1-i}$ satisfy\label{la_1-i}
%$
%1\leq\la_{1-i}\leq p-1
%$
%\end{ass}
%We refer you to \cite{I} for notation.
%S. Ullom proves that if Assumption~\ref{ass:Ullomcond} holds then, for even $i$,
%\begin{equation}\label{eq:Ullom}
%e_i V_n\cong \frac{\Z}{p^n \Z}\oplus (\frac{\Z}{p^{n-1} \Z})^{\la_{1-i}-1}.
%\end{equation}
%This yields, under the same assumption, that
%\begin{equation}\label{eq:Ullom2}
%V_n^+\cong (\frac{\Z}{p^n \Z})^{r(p)} \oplus
%%(\frac{\Z}{p^{n-1}\Z})^{\sum_{even i}(\la_{1-i}-1)}.
%(\frac{\Z}{p^{n-1}\Z})^{\la-r(p)},
%\end{equation}
%where
%\[
%\la=\sum_{i=1,\, i \text{ even}}^{r(p)} \la_{1-i}
%\]
%Hence,
%%under the assumption $\la_{1-i}=0$ or $1$
%when $\la=r$
%we get the first of Kervaire and Murthy's conjectures.
%Note however, that if $\la>r$ the conjecture
%%~\ref{eq:conj2}
%is false.

In this paper we use a different approach.
Instead of directly studying $\Z C_{p^n}$ we
study
\[
A_n:=\frac{\Z[x]}{\Big(\frac{x^{p^n}-1}{x-1}\Big)}
\]
One can prove that $\pic \Z C_{p^n}\cong \pic A_n$. We then construct a pull back similar to Kervaire and Murthy's, but
based on the ring $A_n$ instead of $\Z C_{p^n}$. This gives us a different description of the group $V_n$. We then
construct an injection $ \alpha:\Z[\z_{n-1}]^* \rightarrow A_n^*$ and we show Kervaire and Murthy's $\V_n$ is isomorphic
to the group
\[
\frac{(A_n /pA_n)^*}{\alpha(\Z[\z_{n-1}]^*) \m p}.
\]
In some sense this may be more natural since $\V_n^+$ which we want to calculate is conjectured to be isomorphic to
$\cl^{(p)} \Q (\z_{n-1})$, not $\cl^{(p)} \Q (\z_{n})$. We proceed by constructing a surjection $\V_n^+\rightarrow\V_{n-1}^+$
and an injection  $\V_{n-1}^+\rightarrow\V_{n}^+$ and this allows us to show that
\[
\V_n^+ \cong \big(\frac{\Z}{p^{n}\Z}\big)^{r_0}\oplus
\big(\frac{\Z}{p^{n-1}\Z}\big)^{r_1 - r_0}\oplus \hdots \oplus
\big(\frac{\Z}{p\Z}\big)^{r_{n-1} - r_{n-2}}
\]
for all semi-regular primes.
The sequence $\{r_k\}$ which we describe in section~\ref{chap:prel} is related to the order of certain groups of units
in $\Z[\z_k]$.

In section~\ref{chap:weak} we then go on and prove a weak version of Kervaire and Murthys conjecture \ref{eq:conj3}.
For a group
$A$, let $A(p):=\{a\in A\: : \: a^p=1\}$. We prove that for semi-regular primes,
$\cl \Q(\z_{n-1})(p)\cong \V_n^+/(\V_n^+)^p$ by
using Kervaire and Murthys injection $\cha \V_n^+ \rightarrow \cl^{(p)} \Q(\z_{n-1})$ and constructing a new injection
$\cl \Q(\z_{n-1})(p)\rightarrow \cha \V_n^+$.

In section 4 we then use this weak version and results from section~\ref{chap:prel} to give a proof of
Kervaire and Murthys conjecture \ref{eq:conj3}. The proof relies on class field theory. Since we already have described
the group $\V_n^+$ we of course also get a description of $\cl^{(p)}\Q(\z_{n-1})$.

In section~\ref{chap:appl} we give some applications of our results. When the Iwasawa invariant $\la$ equals $r(p)$,
the index of regularity, we show that all three of Kervaire and Murthy's conjectures hold and that
the Iwasawa invariant
$\nu$ also equals $r(p)$. In fact, $\la=r(p)$ is equivalent to that all three conjectures hold.

We then move on and consider the same assumption Ullom used in \cite{U2}. We show that this leads to that
$\la=r_1$ and prove that $V_n^+=\V_n^+$ in this case too.

Finally we use our results to give some information about the structure of the group of units in $\Z[\z_n]$.

\section{Preliminaries}\label{chap:prel}
We start this section by defining some rings that in some sense are close to $\Z C_{p^n}$. We discuss why
we can and want to work with these rings instead of $\Z C_{p^n}$ and go on get an exact Mayer-Vietoris
sequence from a certain pullback of these rings.

Let for $k \geq 0$ and  $l\geq 1$ \label{A_k,l}
\[
A_{k,l}:=\frac{\Z[x]}{\big( {\frac{x^{p^{k+l}}-1}{x^{p^k} -1}} \big) }
\]
and \label{D_k,l}
\[
D_{k,l}:=A_{k,l} \m p.
\]
We denote the class of $x$ in $A_{k,l}$ by $x_{k,l}$ and in $D_{k,l}$ by $\bar{x}_{k,l}$. Sometimes
we will, by abuse of notation, just denote classes by $x$.
Note that $A_{n,1}\cong \Z[\z_{n}]$ and that
\[
D_{k,l} \cong \frac{\F_p[x]}{(x-1)^{p^{k+l}-p^k}}.
\]

By a generalization of Rim's theorem
(see for example \cite{ST94}) $\pic \Z C_{p^n}\cong \pic A_{0,n}$ for all
$n \geq 1$ so for our purposes we can just as well work with $A_{0,n}$ instead of directly with
$\Z C_{p^n}$. It is easy to see that there exists a pullback diagram
\begin{equation}\label{pullbackdiagram}
\xymatrix@=40pt{
A_{k,l+1} \ar[r]^{i_{k,l+1}} \ar[d]_{j_{k,l+1}}  & \Z[\z_{k+l}]
\ar[dl]^{N_{k,l}}
\ar[d]^{f_{k,l}}\\
A_{k,l} \ar[r]^{g_{k,l}}         &  D_{k,l}
}
\end{equation}
where
$i_{k,l+1}(x_{k,l+1})=\z_{k+l}$,
$j_{k,l+1}(x_{k,l+1})=x_{k,l}$,
$f_{k,l}(\z_{k+l})=\bar{x}_{k,l}$ and
$g_{k,l}$ is just taking classes modulo $p$. The norm-maps $N_{k,l}$  will
be constructed later in this paper. These maps are really the
key to our methods.

The pullback~\ref{pullbackdiagram}
induces a Mayer-Vietoris exact sequence
\begin{multline*}%\label{eq:MV1}
\Z[\z_{n}]^{*}\oplus A_{0,n}^{*}
%\stackrel{\beta}{\rightarrow}
\rightarrow
D_{0,n}^{*} \rightarrow \pic A_{0,n+1} %\rightarrow\\
\rightarrow\pic \Z[\z_{n}]\oplus \pic A_{0,n}\rightarrow
\pic D_{0,n},
\end{multline*}

Since $D_{0,n}$ is local, $\pic D_{0,n}=0$ and since $\Z[\z_{n}]$ is a
Dedekind ring,
$\pic \Z[\z_{n}] \cong \cl \Q (\z_{n})$.
By letting $V_n$ be the cokernel
\[
V_n:=\frac{D_{0,n}^{*}}{\Image\{\Z[\z_{n}]^{*}\times A_{0,n}^{*}\rightarrow D_{0,n}^{*}\}}
\]
we get an exact sequence
\[
0\rightarrow V_n
\rightarrow \pic A_{0,n+1} \rightarrow
\cl \Q (\z_{n})\oplus \pic A_{0,n}\rightarrow 0.
\]
Note that definition of $V_n$ is slightly different from the one from
\cite{K-M}.
By abuse of notation, let $c$ denote the automorphisms on
$A_{k,l}^*$, $\Z[\z_n]^*$ and $D_{k,l}^*$ induced by
$c(t)=t^{-1}$ for $t=x_{k,l}$, $t=\z_n$ and $t=\bar{x}_{k,l}$ respectively. We
also denote the maps induced on $\V_n$ and $V_n$ by $c$.  \label{V_n^+-2}

Before moving on we need to introduce the map $N_{k,l}$.
An element $a \in A_{k,l+1}$ can be uniquely represented as a pair
$(a_l,b_l) \in \Z[\z_{k+l}] \times A_{k,l}$. Using a similar argument
on $b_l$, and then repeating this, we find that $a$ can also be uniquely
represented as an $(l+1)$-tuple $(a_l,\ldots,a_m, \ldots,a_0)$ where
$a_m \in \Z[\z_{k+m}]$. In the rest of this paper we will identify an element
of $A_{k,l+1}$ with both its representations as a pair or an $(l+1)$-tuple.

For $k\geq 0$ and $l \geq 1$ let
$\tilde{N}_{k+l,l}:\Z[\z_{k+l}]\rightarrow \Z[\z_k]$ denote the usual norm.

\begin{prop}\label{prop:Norms}
For each $k \geq 0$ and $l \geq 1$ there exists a multiplicative map
$N_{k,l}$ such that the diagram
\[
\xymatrix@=40pt{
  &  \Z[\z_{k+l}] \ar[d]^{f_{k,l}} \ar[dl]_{N_{k,l}}\\
A_{k,l} \ar[r]^{g_{k,l}}         &  D_{k,l}
}
\]
is commutative.
Moreover, if $a\in \Z[\z_{k+l}]$, then
\[
N_{k,l}(a)=(\tilde{N}_{k+l,1}(a),N_{k,l-1}(\tilde{N}_{k+l,1}(a)))=
(\tilde{N}_{k+l,1}(a),\tilde{N}_{k+l,2}(a),\hdots,\tilde{N}_{k+l,l}(a)).
\]
\end{prop}
The construction of $N_{k,l}$
can be found in \cite{ST97}. Since it may not be well known we will  for completeness
repeat it here.
Before this we notice an immediate
consequence of the commutativity of the diagram in Proposition~\ref{prop:Norms}.
\begin{cor}
$V_n=\frac{D_{0,n}^{*}}{\Image\{A_{0,n}^{*}\rightarrow D_{0,n}^{*}\}}$
\end{cor}
\def\proofname{\bf Proof of Proposition 2.1.}
\begin{proof}
The maps $N_{k,l}$ will be constructed inductively. If $i=1$ and $k$
is arbitrary, we have  $A_{k,1}\cong \Z[\z_{k}]$ and we define
$N_{k,1}$ as the usual norm map $\tilde{N}_{k+1,1}$. Since
$\tilde{N}_{k+1,1}(\z_{k+1})=\z_k$ we only need to prove that our
map is additive modulo $p$, which follows from the lemma below.
\begin{lemma}
For $k \geq 0$ and $l\geq 1$ we have
\begin{itemize}\vspace{-11pt}
\item [{\it i}$)$]
          $A_{k+1,l}$ is a free $A_{k,l}$-module under
          $x_{k,l}\mapsto x_{k+1,l}^p$.
\item [{\it ii}$)$]
          The norm map $N:A_{k+1,l}\rightarrow A_{k,l}$, defined by
          taking the determinant of the multiplication operator, is additive
          modulo $p$.
\end{itemize}
\end{lemma}
This is Lemma 2.1 and Lemma 2.2 in \cite{ST97} and proofs can be found there.

Now suppose $N_{k,j}$ is constructed for all $k$ and all $j\leq l-1$. Let
$\varphi=\varphi_{k+1,l}:\Z[\z_{k+l}]\rightarrow A_{k+1,l}$ be defined by
$\varphi(a)=(a,N_{k+1,l-1}(a))$. It is clear that $\varphi$ is multiplicative.
From the lemma above we have a norm map $N:A_{k+1,l}\rightarrow A_{k,l}$.
Define $N_{k,l}:=N \circ \varphi$. It is clear that $N_{k,l}$ is
multiplicative. Moreover,
$N_{k,l}(\z_{k+l})=N(\z_{k+l},x_{k+1,l-1})=N(x_{k+1,l})=x_{k,l}$, where
the latter equality follows by a direct computation. To prove that our map
makes the diagram in the proposition above commute, we now only need
to prove it is additive modulo $p$. This also follows by a direct
calculation once you notice that
\[
\varphi(a+b)-\varphi(a)-\varphi(b) =
\frac{x_{k+1,l}^{p^{k+l+1}}-1}{x_{k+1,l}^{p^{k+l}}-1} \cdot r,
\]
for some $r \in A_{k+1,l}$.

Regarding the other two equalities in Proposition~\ref{prop:Norms},
it is clear that the second one follows from the first. The first equality
will follow from the lemma below.

\begin{lemma}
The diagram
\[
\xymatrix@=40pt{
\Z[\z_{k+l}] \ar[r]^{{\tilde{N}}_{k+l,1}} \ar[d]_{N_{k,l}} &
         \Z[\z_{k+l-1}] \ar[d]_{N_{k-1,l}} \\
A_{k,l} \ar[r]^{N} & A_{k-1,l}\\
}
\]
is commutative
\end{lemma}
\begin{proof}
Recall that the maps denoted $N$ (without subscript) are the usual norms
defined by the determinant of the multiplication map.
An element in $A_{k,l}$
can be represented as a pair $(a,b)\in \Z[\z_{k+l-1}]\times A_{k,l-1}$ and
an element in  $A_{k-1,l}$
can be represented as a pair $(c,d)\in \Z[\z_{k+l-2}]\times A_{k-1,l-1}$.
If $(a,b)$ represents an element in $A_{k,l}$ one can,
directly from the definition, show that
$N(a,b)=({\tilde{N}}_{k+l-1,1}(a),N(b)) \in A_{k-1,l}$.
We now use induction on $l$.
If $l=1$ the statement is well known. Suppose the
diagram corresponding to the one above, but with $l$ replaced
by $l-1$, is commutative for all $k$.
If $a\in \Z[\z_{k+l}]$ we have
\[
N(N_{k,l}(a))=N(N((a,N_{k+1,l-1}(a)))=({\tilde{N}}_{k+l,2}(a),N(N(N_{k+1,l-1}(a))))
\]
and
\[
N_{k-1,l}(N(a))=( {\tilde{N}}_{k+l,2}(a) ,N(N_{k,l-1}({\tilde{N}}_{k+l,1}(a)))).
\]
By the induction hypothesis applied to $(k+1, l-1)$ we get
$N_{k,l-1}\circ {\tilde{N}}_{k+l,1}=N\circ N_{k+1,l-1}$ and this proves the lemma.
\end{proof}
With the proof of this Lemma the proof of Proposition~\ref{prop:Norms} is complete.
\end{proof}

We will now use our the maps $N_{k,l}$ to get an inclusion
of $\Z[\z_{k+l-1}]^*$ into $A_{k,l}^*$.
Define $\varphi_{k,l}:\Z[\z_{k+l-1}]^*\rightarrow A_{k,l}^*$ be the injective
group homomorphism defined by
$\e \mapsto (\e,N_{k,l}(\e))$. By
Proposition~\ref{prop:Norms},
$\varphi_{k,l}$ is well defined. For future use we record this in a lemma.
\begin{lemma}\label{lemma:A*decomp}
Let $B_{k,l}$ be the subgroup of $A_{k,l}^*$ consisting of elements
$(1,b)$, $b \in A_{k,l-1}^*$.
Then $A_{k,l}^*\cong \Z[\z_{k+l-1}]^* \times B_{k,l}$
\end{lemma}
In what follows, we identify
$\Z[\z_{k+l-1}]^*$ with its image in $A_{k,l}^*$.

Before we move on we will state a technical lemma which is Theorem I.2.7 in \cite{ST98}.
\begin{lemma}\label{lemma:gker}
Let $a\in \Z[\z_{k+l-1}]$. Then
$g_{k,l}(a, N_{k,l-1}(a))=1\in D_{k,l}$ if
and only if $a \equiv 1 \m \lambda_{k+l-1}^{p^{k+l}-p^k}\}$.
In particular,
\[
\kernel ({g_{k,l}}{|_{\Z[\z_{k+l-1}]^*}}) =
\{\e \in \Z[\z_{k+l-1}]^*\, : \, \e \equiv 1 \m \lambda_{k+l-1}^{p^{k+l}-p^k}\}
\]
\end{lemma}
We will not repeat the proof here, but since the technique used is
interesting we will indicate the main idea. If $a\in \Z[\z_{k+l-1}]^{*}$
and $g_{k,l}(a)=1$ we get that $a\equiv 1 \m p$ in $\Z[\z_{k+l-1}]$,
 $N_{k,l-1}(a)\equiv 1 \m p$ in $A_{k,l-1}$ and that
$f_{k,l-1}\big(\frac{a-1}{p}\big) =
g_{k,l-1}\big(\frac{N_{k,l-1}(a)-1}{p}\big)$. Since the norm map
commutes with $f$ and $g$ this means that
$N_{k,l-1}(\frac{a-1}{p}) \equiv \frac{N_{k,l-1}(a)-1}{p}$. The latter
is a congruence in $A_{k,l-1}$ and by the same method as above we
deduce a congruence in
$\Z[\z_{k+l-2}]$ and a
 congruence in $A_{k,l-2}$. This can be repeated $l-1$ times
until we get a congruence in $A_{k,1}\cong \Z[\z_k]$. The last
congruence in general looks pretty complex, but can be analyzed
and gives us the neccesary information.

If for example $l=2$, we get after just one step  $a \equiv 1 \m p$ in
$\Z[\z_{k+1}]$, $N(a)\equiv 1 \m p$ and
$N(\frac{a-1}{p}) \equiv \frac{N(a)-1}{p} \m p$ in
$A_{k,1}\cong \Z[\z_k]$, where $N$ is the usual norm. By viewing $N$ as
a product of automorphisms, recalling that $N$ is additive modulo $p$ and that
the usual trace of any element of $\Z[\z_{k+1}]$ is divisible by $p$, we get
that $N(a) \equiv 1 \m p^2$ and hence that
$N(\frac{a-1}{p}) \equiv 0 \m p$. By analyzing how the norm acts one
can show that this means that $a \equiv 1 \m \la_{k+1}^{p^{k+2}-p^k}$

In the rest of this paper we paper will only need the the rings $A_{k,l}$ and $D_{k,l}$ in the case
$k=0$. Therefore we will simplify the notation a little by setting $A_l:=A_{0,l}$,
$D_l:=D_{0,l}$, $g_l:=g_{0,l}$, $f_l:=f_{0,l}$, $i_l:=i_{0,l}$, $j_l:=j_{0,l}$
and $N_l:=N_{0,l}$.
%%%%%%%%%%%New indexes starts from here:

Now define $\V_n$ as
\[
\V_n:=\frac{\tilde{D}_{n}^{*}}{\Image\{\tilde{\Z}[\z_{n-1}]^*\rightarrow \tilde{D}_{n}^{*}\}},
\]
where $\tilde{\Z}[\z_{n-1}]^*$ are the group of all units $\e$ such that
$\e \equiv 1 \m \la_{n-1}$, where $\la_n$ denotes the ideal $(\z_n -1)$,
and $\tilde{D}_{n}^{*}$ are the units that are congruent to $1$ modulo
the class of $(\bar{x}-1)$ in $D_{n}^{*}$.
This definition is equivalent to the definition in \cite{K-M} by the following Proposition.
%by Lemma~\ref{lemma:Nsurj},
%$N:\Z[\z_{n}]^*\rightarrow \Z[\z_{n-1}]^*$ is surjective when $p$ i semi-regular.
\begin{prop}\label{Vn=Vn}
The two definitions of $\V_n$  coincide.
\end{prop}
\begin{proof}
The kernel of the surjection $(\F_p[x]/(x-1)^{p^n})^*\rightarrow
(\F_p[x]/(x-1)^{p^n-1})^*=D_n^*$ consists of units congruent to 1
mod $(x-1)^{p^n-1}$.
Let $\eta:=\z_{n}^{\frac{p^{n+1}+1}{2}}$. Then $\eta^2=\z_{n}$ and $c(\eta)=\eta^{-1}$.
Let
$
\e:=\frac{\eta^{p^{n}+1}-\eta^{-(p^{n}+1)}}{\eta-\eta^{-1}}.
$
One can by a direct calculation show that
%see lic.tex for details
$\e=1+(\z_{n}-1)^{p^{n}-1}+t(\z_{n}-1)^{p^{n}}$
for some $t\in \Z[\z_{n}]$.
If $a=1+a_{p^{n}-1}(x_{n}-1)^{p^{n}-1}\in (\F_p[x]/(x-1)^{p^n})^*$,
$a_{p^{n}-1}\in \F_p^*$,
Then it is just a matter of calculations to show that
$a=f_{n}(\e)^{a_{p^{n}-1}}$. This shows that
 $(\F_p[x]/(x-1)^{p^n})^*/f_n'(\Z[\z_{n}]^*)\cong
(\F_p[x]/(x-1)^{p^n-1})^*/f_n(\Z[\z_{n}]^*)$.
Since
\[
\xymatrix@=40pt{
& \Z[\z_{n}]^*  \ar[dl]_{N}
\ar[d]^{f}\\
\Z[\z_{n-1}]^*  \ar[r]^{g}         &    \tilde{D}_{n}^{*+}
}
\]
is commutative and $N$ (which is the restriction of the usual norm-map) surjective
when $p$ is semi-regular (Lemma~\ref{lemma:Nsurj})
the proposition follows.
\end{proof}

Let $\V_n^+:=\{v\in \V_n\; :\; c(v)=v\}$. What we want to do is to find the structure of $\V_n^+$.
For $n \geq 0$ and $k\geq 0$, define
\[
U_{n,k}:=\{ real\,\,\e \in \Z[\z_n]^* : \e \equiv 1 \m \la_n^{k}\}.
\]
One of our main results is the
following proposition.
\begin{prop}\label{prop:Vn-orderind}
If $p$ is semi-regular,
$|\V_n^+|=|\V_{n-1}^+| \cdot |U_{n-1,p^{n}-1} / (U_{n-1,p^{n-1}+1})^{p}|$.
\end{prop}
Here $U^{p}$ denotes the group of $p$-th powers of elements of the group $U$. Note that we in this paper sometimes
use the notation $R^n$ for $n$ copies of the ring (or group) $R$. The context will make it clear which one of these two
things we mean.

For $k=0,1,\hdots$, define $r_k$ by
\[
|U_{k,p^{k+1}-1} / (U_{k,p^{k}+1})^{p}|=p^{r_k}.
\]
By Lemma 2 in \cite{ST94} we get that $U_{k,p^{k+1}-1}=U_{k,p^{k+1}}$ and since the
the $\la_n$-adic valuation of $\e -1$, where $\e$ is a real unit, is even,
$U_{k,p^{k+1}}=U_{k,p^{k+1}+1}$. We hence have
\begin{lemma}\label{lemma:p^k-1=p^k+1}
$U_{k,p^{k+1}-1}=U_{k,p^{k+1}+1}$.
\end{lemma}

One can prove that $r_0=r(p)$, the index of irregularity, since if the $\la_0$-adic valuation
of $\e\in \Z[\z_0]^{*+}$ is less than $p-1$, then local considerations show that the extension
$\Q(\z_0)\subseteq\Q(\z_0,\sqrt[p]{\e})$ is ramified. The result then follows from the fact
that
\[
\frac{U_{0,p-1}}{(U_{0,2})^p}\cong\frac{S_0}{pS_0}
\]
where $S_0$ is the $p$-class group of $\Q(\z_0)$.

Before the proof of  Proposition~\ref{prop:Vn-orderind} we will state and a lemma, which
is well-known.

\begin{lemma}\label{lemma:Nsurj}
If $p$ is semi-regular  $N_{n-1}:\Z[\z_{n-1}]\rightarrow A_{n-1}$ maps
$U_{n-1,1}$ surjectively onto $U_{n-2,1}$.
\end{lemma}
\def\proofname{\bf Proof.}

\def\proofname{\bf Proof of Proposition~\ref{prop:Vn-orderind}.}
\begin{proof}
\def\proofname{\bf Proof.}
In a similar way as the ideal $\la_n:=(\z_n-1)$ equal the ideal $(\z_n-\z_n^{-1})$ in $\Z[\z_n]$
one can show that that $(\bar{x}-1)=(\bar{x}-\bar{x}^{-1})$ in $D_{n}$. It is easy to show that
$\tilde{D}_{n}^{*+}$ can be represented by elements
$1+a_2(\bar{x}-\bar{x}^{-1})^2+a_4(\bar{x}-\bar{x}^{-1})^4+\hdots + a_{p^{n}-3}(x-x^{-1})^{p^{n}-3}$,
$a_i\in \F_p$. Hence $|\tilde{D}_{n}^{*+}|=p^{(p^{n}-3)/2}$. We want to evaluate
\[
|\tilde{D}_{n}^{*+}|/|g_{n}(U_{n-1,1})|.
\]
By Lemma~\ref{lemma:gker} we have
\[
g_{n}(U_{n-1,1} ) \cong \frac{U_{n-1,1} }{U_{n-1,p^{n}-1} }.
\]
Since $g_{n}(U_{n-1,1} ) \subseteq
g_{n}(\Z[\z_{n-1}]^{*+}) \subseteq
\tilde{D}_{n}^{*+}$ the group
$U_{n-1,1} /U_{n-1,p^{n}-1} $ is finite.
Similarly $\Z[\z_{n-1}]^{*+}/U_{n-1,p^{n}-1} $ is finite.
This shows that $\Z[\z_{n-1}]^{*+}/U_{n-1,1} $ is finite since
\[
\big|\frac{\Z[\z_{n-1}]^{*+}}{U_{n-1,1} }\big|
\big|\frac{U_{n-1,1} }{U_{n-1,p^{n}-1} }\big|=
\big|\frac{\Z[\z_{n-1}]^{*+}}{U_{n-1,p^{n}-1} }\big|.
\]

We can write
%\begin{eqnarray*}\label{eq:U/U}
\begin{equation}\label{eq:U/U}
\begin{split}
\big|\frac{U_{n-1,1} }{U_{n-1,p^{n}-1} }\big|
=
%&=&
\big|\frac{U_{n-1,1} }{U_{n-1,p^{n-1}-1} }\big|
\big|\frac{U_{n-1,p^{n-1}-1} }{U_{n-1,p^{n-1}+1} }\big|
\big|\frac{U_{n-1,p^{n-1}+1} }{U_{n-1,p^{n}-1} }\big|=\\
=
%&=&
\big|\frac{U_{n-1,1} }{U_{n-1,p^{n-1}-1} }\big|
\big|\frac{U_{n-1,p^{n-1}-1} }{U_{n-1,p^{n-1}+1} }\big|
\big|\frac{U_{n-1,p^{n-1}+1} /(U_{n-1,p^{n-1}+1} )^p}{U_{n-1,p^{n}-1} /(U_{n-1,p^{n-1}+1} )^p}\big|=\\
=
%&=&
\big|\frac{U_{n-1,1} }{U_{n-1,p^{n-1}-1} }\big|
\big|\frac{U_{n-1,p^{n-1}-1} }{U_{n-1,p^{n-1}+1} }\big|
\big|\frac{U_{n-1,p^{n-1}+1} }{(U_{n-1,p^{n-1}+1} )^p}\big|
\big|\frac{U_{n-1,p^{n}-1} }{(U_{n-1,p^{n-1}+1} )^p}\big|^{-1}
\end{split}
\end{equation}
%\end{eqnarray*}
%
%
By Dirichlet's theorem on units we have
$(\Z[\z_{n-1}]^*)^+\cong \Z^{\frac{p^{n}-p^{n-1}}{2}-1}$
Since all quotient groups involved are finite we get that
$U_{n-1,1}$, $U_{n-1,p^{n}-1}$, $U_{n-1,p^{n-1}-1}$
and $U_{n-1,p^{n-1}+1}$ are all isomorphic to
$\Z^{\frac{p^{n}-p^{n-1}}{2}-1}$.
The rest of the proof is devoted to the analysis of the four right
hand factors of~\ref{eq:U/U}.

Obviously,
\[
\frac{U_{n-1,p^{n-1}+1} }{(U_{n-1,p^{n-1}+1} )^p}\cong
\frac{\Z^{\frac{p^{n}-p^{n-1}}{2}-1}}{(p\Z)^{\frac{p^{n}-p^{n-1}}{2}-1}}
\cong C_p^{\frac{p^{n}-p^{n-1}}{2}-1}.
\]
This shows that
\[
\big|\frac{U_{n-1,p^{n-1}+1} }{(U_{n-1,p^{n-1}+1} )^p}\big|=
p^{\frac{p^{n}-p^{n-1}}{2}-1}.
\]

We now turn to the second factor of the right hand side of~\ref{eq:U/U}.
We will show that this number is $p$ by finding a unit
$\e \not\in {U_{p^{n-1}+1} }$ such that
\[
<\e>=\frac{U_{n-1,p^{n-1}-1} }{U_{n-1,p^{n-1}+1} }.
\]
%\[
%<\e>=\frac{U_{n-1,p^{n-1}-1} }{U_{n-1,p^{n-1}+1} }.
%\]
Since the $p$-th power of any unit in $U_{n-1,p^{n-1}-1} $ belongs to
$U_{n-1,p^{n-1}+1} $ this is enough. Let $\z=\z_{n-1}$ and
$\eta:=\z^{\frac{p^{n}+1}{2}}$. Then $\eta^2=\z$ and $c(\eta)=\eta^{-1}$.
Let
$
\e:=\frac{\eta^{p^{n-1}+1}-\eta^{-(p^{n-1}+1)}}{\eta-\eta^{-1}}.
$
Then $c(\e)=\e$ and one can by direct calculations show that
$\e$ is the unit we are looking for.
%
%ALLA BERÄKNINGAR FINNS I lic.tex
%
%

We now want to calculate
\[
\big|\frac{U_{n-1,1} }{U_{n-1,p^{n-1}-1} }\big|.
\]
Consider the
commutative diagram
\[
\xymatrix@=40pt{
& \Z[\z_{n-1}]^* \ar[dl]_{N_{n-1}}
\ar[d]^{f_{n-1}}\\
A_{n-1}^* \ar[r]^{g_{n-1}}         &  D_{n-1}^*
}
\]
It is clear that
$f_{n-1}(U_{n-1,1} )\subseteq \tilde{D}_{n-1}^{*+}$ and that
$g_{n-2}(U_{n-2,1} )\subseteq \tilde{D}_{n-1}^{*+}$.
By Lemma~\ref{lemma:Nsurj}
we have a commutative diagram
\[
\xymatrix@=40pt{
& U_{n-1,1}  \ar[dl]_{N}
\ar[d]^{f}\\
U_{n-2,1}  \ar[r]^{g}         &    \tilde{D}_{n-1}^{*+}
}
\]
where $N$ is surjective. Clearly, $f(U_{n-1,1} )=g(U_{n-2,1} )$.

It is easy to see that
$\kernel (f)=U_{n-1,p^{n-1}-1} $ so by above
\[
\frac{U_{n-1,1} }{U_{n-1,p^{n-1}-1} } \cong
   f(U_{n-1,1} ) = g(U_{n-2,1} ).
\]
Now recall that by definition
$\V_{n-1}^+=\tilde{D}_{n-1}^{*+}/g(U_{n-2,1} )$. Hence
\[
\big|\frac{U_{n-1,1} }{U_{n-1,p^{n-1}-1} }\big|=
    |g(U_{n-2,1} )| = |\tilde{D}_{n-1}^{*+}||\V_{n-1}^+|^{-1}=
p^{\frac{p^{n-1}-3}{2}}|\V_{n-1}^+|^{-1}.
\]
This finally gives
\begin{eqnarray*}
|\V_{n}^+|&=&|\tilde{D}_{n}^{*+}||g(U_{n-1,1} )|^{-1}=\\
&=&
p^{\frac{p^{n}-3}{2}} \cdot
p^{-\frac{p^{n-1}-3}{2}} \cdot |\V_{n-1}^+| \cdot
p^{-1} \cdot
p^{-\frac{p^{n}-p^{n-1}}{2}+1} \cdot
\big|\frac{U_{n-1,p^{n}-1} }{(U_{n-1,p^{n-1}+1} )^p}\big|=\\
&=&
|\V_{n-1}^+| \cdot \big|\frac{U_{n-1,p^{n}-1} }{(U_{n-1,p^{n-1}+1} )^p}\big|
\end{eqnarray*}
which is what we wanted to show.

\end{proof}
\def\proofname{\bf Proof.}

\begin{prop} The sequence
$\{r_k\}$ is non-decreasing, bounded by
the Iwasawa invariant $\la$ and
$|\V_n^+|=p^{r_0+ r_1 +\hdots + r_{n-1}}$.
\end{prop}
\begin{proof}
Recall that $(\la_{k})=(\la_{k+1}^p)$ as ideals in $\Z[\z_{k+1}]$.
By Lemma~\ref{lemma:p^k-1=p^k+1}, the inclusion of $\Z[\z_k]$ in $\Z[\z_{k+1}]$ induces an inclusion
of $U_{k,p^{k+1}-1} =U_{k,p^{k+1}+1} $ into $U_{k+1,p^{k+2}+p} \subseteq U_{k+1,p^{k+2}-1} $.
Since a $p$-th power in
$\Z[\z_k]$ obviously is a $p$-th power in $\Z[\z_{k+1}]$ we get an homomorphism of
\begin{equation}\label{eq:U-homo}
\frac{U_{k,p^{k+1}-1} }{(U_{k,p^{k}+1} )^{p}}
\rightarrow \frac{U_{k+1,p^{k+2}-1} }{(U_{k+1,p^{k+1}+1} )^{p}}.
\end{equation}
If $\e\in U_{k,p^{k+1}-1} $ is a not $p$-th power in $\Z[\z_{k}]$ then one can show that
$\Q(\z_{k})\subseteq \Q(\z_{k}, \e)$ is an unramified extension of degree $p$. If $\e$ would
be a $p$-th power in $\Z[\z_{k+1}]$ we would get $\Q(\z_{k+1})=\Q(\z_{k}, \e)$ which is impossible
since $\Q(\z_{k})\subseteq Q(\z_{k+1})$ is ramified. Hence the homomorphism~\ref{eq:U-homo}
%\[
%\frac{U_{k,p^{k+1}-1} }{(U_{k,p^{k}+1} )^{p}}
%\rightarrow \frac{U_{k+1,p^{k+2}-1} }{(U_{k+1,p^{k+1}+1} )^{p}}
%\]
is injective.
This shows that the sequence
$\{r_k\}$ is non-decreasing.

%Let $\cl^{p}(\Q(\z_n))$ be the $p$-Sylow subgroup of the class group of $\Q(\z_n)$ and
%let $r_n=r_n(p)$ be the $p$-rank of $\cl^{p}(\Q(\z_n))$. It is known that $r_0(p)=r(p)$,
%the index of regularity of the semi-regular prime $p$.
%It turns out that
%with the help of some tricks, class field theory and Kummers pairing one can prove the
%following lemma.
%\begin{lemma}\label{lemma:U/Up-order}
%$|U_{n-1,p^{n}-1}  / (U_{n-1,p^{n-1}+1} )^{p}|=p^{r_{n-1}}$
%\end{lemma}

Since it is known by for example \cite{K-M} that $|\V_1^+|=p^{r_0}$,
by induction and Proposition~\ref{prop:Vn-orderind} we now immediately get:
$|\V_n^+|=p^{r_0+ r_1 +\hdots + r_{n-1}}$.

Assume now that $r_N >\la$ for some $N$. Then it follows that $|\V_n^+|$ grows faster
than $p^{\la n}$ and this contradicts to that of  $|\V_n^+|< |\cl^{(p)} \Q (\z_{n})|=p^{\la n +\nu}$.
This observation completes the proof.
\end{proof}

%Consider the surjection $\pi_n:\V_n^+\rightarrow \V_{n-1}^+$ from ~\ref{prop:Vsurj}.
%Now consider the following lemma.

We now go on and prove the following proposition:
\begin{prop}\label{prop:Vsurj}
There exists a surjection
$\pi_n:\V_n^+ \rightarrow \V_{n-1}^+$.
\end{prop}
\begin{proof}
The canonical surjection
$j_{n}:A_{n}\rightarrow A_{n-1}$ can be considered mod $(p)$ and hence
yields a surjection $\bar{j}_{n}:D_{n} \rightarrow D_{n-1}$.
Suppose that $\bar{u}\in D_{n-1}^{*+}$, $\bar{v}\in D_{n}^{*+}$,
$\bar{j}_{n}(\bar{v})=\bar{u}$ and that $\bar{v}=g_{n}(v)$, where
$v=(\e,N_{n-1}(\e))$, $\e\in \Z[\z_{n-1}]$. Then $j_{n}(v)=N_{n-1}(\e)$,
and $\bar{u}=\bar{j}_{n}(\bar{v})=\bar{j}_{n}g_{n}N_{n-1}(\e)$.
But $N_{n-1}(\e)=(\tilde{N}_{n-1,1}(\e),N_{n-2}\tilde{N}_{n-1,1}(\e))$ by
Proposition~\ref{prop:Norms}. In other words, if $\bar{v}$ represents $1$ in
$\V_n$, then $\bar{j}_{n}(\bar{v})$ represents $1$ in $\V_{n-1}$ so the map
$\bar{j}_{n}$ induces a well defined surjection $\V_n^+\rightarrow \V_{n-1}^+$.
%%%%%%%%%%%
%
%
%Some kind of injectivity of \bar{j} restricted to g(\Z[\z_{n}]) is below. Can be useful.
%
%
%
%Suppose that under
%this latter surjection, $u \in D_{n-1}^{*+}$ is the
%image of $v \in D_{n}^{*+}$ and suppose $u=g_{n-1}((\e,N_{n-2}(\e)))$
%for some $\e \in U_{n-2,1}  \subset \Z[z_{n-2}]$.
%For some $a \in A_{n}$, $v=g_{n}(a)$ and
%$(\e,N_{n-2}(\e))=j_{n}(a)$. Since $p$ is semi-regular we know from
%Lemma~\ref{lemma:Nsurj}
%that the norm map $N_{n-1}$ resticted to
%$U_{n-1,1} $ is surjective onto
%$U_{n-2,1} $ and acts as the usual norm
%$\tilde{N}_{n-1,1}$. Hence there exists $\e' \in U_{n-1,1} $ such that
%$N_{n-1}(\e')=(\e,N_{n-2}(\e))$. This means that
%$(\e',N_{n-1}(\e')) \in A_{n}^{*+}$ maps to $(\e,N_{n-2}(\e))$
%under $j_{n}$. Since $f_{n-1}(\e')=g_{n-1}N_{n-1}(\e')=u$ and all the
%maps come from a pullback we get that $a=(\e',N_{n-1}(\e'))$, that is, $v$
%is the image of a unit in $U_{n-1,1} $. This proves the surjection
%$\V_n^+ \rightarrow \V_{n-1}^+$ is well defined.
\end{proof}
\def\proofname{\bf Proof.}

It is now not hard to find the kernal of $\pi_n$.
\begin{prop}\label{prop:pi-kernal}
For any semi-regular prime $p$, $\ker \pi_n \cong (\Z/p\Z)^{r_{n-1}}$.
\end{prop}
\def\proofname{\bf Proof.}
\begin{proof}
Proposition~\ref{prop:Vn-orderind} and the definition of $r_n$ clearly implies that
$|\ker \pi_n |=p^{r_{n-1}}$. We need to prove that any element in $\ker \pi_n$
has order at most $p$.
%Direct proof
Suppose that in the surjection $D_{n}^{*+}\rightarrow D_{n-1}^{*+}$, the element
$u \in D_{n-1}^{*+}$ is the
image of $v \in D_{n}^{*+}$ and suppose $u=g_{n-1}((\e,N_{n-2}(\e)))$
for some $\e \in U_{n-2,1}  \subset \Z[\z_{n-2}]$.
For some $a \in A_{n}$, $v=g_{n}(a)$ and
$(\e,N_{n-2}(\e))=j_{n}(a)$. Since $p$ is semi-regular we know from
Lemma~\ref{lemma:Nsurj}
that the norm map $N_{n-1}$ resticted to
$U_{n-1,1} $ is surjective onto
$U_{n-2,1} $ and acts as the usual norm
$\tilde{N}_{n-1,1}$. Hence there exists $\e' \in U_{n-1,1} $ such that
$N_{n-1}(\e')=(\e,N_{n-2}(\e))$. This means that
$(\e',N_{n-1}(\e')) \in A_{n}^{*+}$ maps to $(\e,N_{n-2}(\e))$
under $j_{n}$. Since $f_{n-1}(\e')=g_{n-1}N_{n-1}(\e')=u$ and all the
maps come from a pullback we get that $a=(\e',N_{n-1}(\e'))$, that is, $v$
is the image of a unit in $U_{n-1,1} $.
Now define $\tilde{D}_{n,(k)}^{*+}:=\{a\in \tilde{D}_{n}^{*+}: a\equiv 1 \m (x-1)^k\}$.
Then
\[
\ker \pi_n =
\frac{\ker\{\tilde{D}_{n}^{*+}\rightarrow \tilde{D}_{n-1}^{*+}\}}{\ker\{\tilde{D}_{n}^{*+}\rightarrow \tilde{D}_{n-1}^{*+}\}\cap g_{n}
(\Z[\z_{n-1}]^{*+})}=
\frac{\tilde{D}_{n,(p^{n-1}-1)}^{*+}}{g_{n}(U_{n-1,p^{n-1}-1}) }.
\]
Now note that if $b \in  \tilde{D}_{n,(p^{n-1})}^{*+}$, then $b^p=1$ so such a unit clearly has order
$p$. We will show that any unit $a\in  \tilde{D}_{n,(p^{n-1}-1)}^{*+}$ can be written as
$a=b g_{n}(\e)^k$ for some $b \in  \tilde{D}_{n,(p^{n-1})}^{*+}$, natural number $k$ and
$\e \in U_{n-1,p^{n-1}-1} $. Then $a^p=b^p g_{n}(\e)^{kp}$ is clearly trivial in
$\ker \pi_n \subseteq  \V_n^+$.
Let $\eta:=\z_{n-1}^{\frac{p^{n}+1}{2}}$. Then $\eta^2=\z_{n-1}$ and $c(\eta)=\eta^{-1}$.
Let
$
\e:=\frac{\eta^{p^{n-1}+1}-\eta^{-(p^{n-1}+1)}}{\eta-\eta^{-1}}.
$
One can by a direct calculation show that
%see lic.tex for details
$\e\in U_{n-1,p^{n-1}-1} \setminus U_{n-1,p^{n-1}+1} $.
In fact, $\e=1+e_{p^{n-1}-1}(\z_{n-1}-\z_{n-1}^{-1})^{p^{n-1}-1}+t(\z_{n-1}-\z_{n-1}^{-1})^{p^{n-1}+1}$
for some non-zero $e_{p^{n-1}-1}\in \Z[\z_{n-2}]$, not divisible by $\la_{n-1}$, and some $t\in \Z$.
Suppose $a=1+a_{p^{n-1}-1}(x_{n-1}-x_{n-1}^{-1})^{p^{n-1}-1}+\ldots \in\tilde{D}_{n,(p^{n-1}-1)}^{*+}$,
$a_{p^{n-1}-1}\in \F_p^*$. Since $e_{p^{n-1}-1}$ is not divisible by $\la_{n-1}$, $g_{n}(\e)\in \F_p^*$
Hence we can choose $k$ such that
$kg_{n}(e_{p^{n-1}-1})\equiv a_{p^{n-1}-1} \m p$.
Then it is just a matter of calculations to show that
$a=b g_{n}(\e)^k$, where $b\in \tilde{D}_{n,(p^{n-1})}^{*+}$, which concludes the proof
\end{proof}

One of our main theorems is the following:
\begin{thm}\label{thm:Vn-semireg}
For every semi-regular prime $p$
\[
\V_n^+ \cong \big(\frac{\Z}{p^{n}\Z}\big)^{r_0}\oplus
\big(\frac{\Z}{p^{n-1}\Z}\big)^{r_1 - r_0}\oplus \hdots \oplus
\big(\frac{\Z}{p\Z}\big)^{r_{n-1} - r_{n-2}}.
\]
\end{thm}

To prove this we need to introduce some techniques
from \cite{K-M}.

Let $P_{0,n}$ be the group of principal fractional ideals in $\Q(\z_{n})$ prime to
$\la_{n}$. Let $H_n$ be the subgroup of fractional ideals congruent to 1 modulo
$\la_{n}^{p^{n}}$. In \cite{K-M}, p. 431, it is proved that there exists a canonical isomorphism
\[
J:\frac{P_{0,n}}{H_n}\rightarrow \frac{(\F_p[x]/(x-1)^{p^n})^{*}}{f_n'(\Z[\z_n]^*)}=:\V_n'.
\]
Now consider the injection $\iota:\Q(\z_{n-1})\rightarrow \Q(\z_{n})$,
$\z_{n-1}\mapsto\z_n^p$. It is clear we get an induced map
$P_{0,n-1}\rightarrow P_{0,n}$. Since $\iota$ map $\la_{n-1}$ to $\la_n^p$ it is easy to see
that we get an induced homomorphism
\[
\alpha'_n:\frac{P_{0,n-1}}{H_{n-1}}\rightarrow \frac{P_{0,n}}{H_n}.
\]
Considered as a map $\alpha'_n:\V_{n-1}'\rightarrow \V_{n}'$ this map acts as
$(\F_p[x]/(x-1)^{p^{n-1}})^{*} \ni x_{n-1} \mapsto x_n^p \in (\F_p[x]/(x-1)^{p^n})^{*}$.
Since $\V_n'\cong \V_n$ (see Proposition~\ref{Vn=Vn}) we can consider this as a homomorphism
$\alpha_n:\V_{n-1}\rightarrow \V_{n}$. Clearly we then get that
$\alpha$ is induced by $x_{n-1}\rightarrow x_{n}^p$
%(x_{n-1}g_{n-1}(\tilde{\Z}[\z_{n-2}]))=x_{n}^p g_{n}(\tilde{\Z}[\z_{n-1}])$.
%, or in other words, $\alpha_n(a_{n-1})=$.
Note however, that $x_{n-1}\mapsto x_{n}^p$ does not induce a homomorphism
$D_{n-1}^*\rightarrow D_n^*$.
\begin{lemma}\label{lemma:Vn-inj}
The map $\alpha_n$ is injective on $\V_{n-1}^+$.
\end{lemma}
\begin{proof}
In this proof, denote $\Q(\z_{n})$ by $F_n$.
Let $L_n$ be the $p$-part of the Hilbert class field of
$F_n$ and let $K_n/F_n$ be the $p$-part  of the ray class field extension
associated with the ray group $H_n$. In other words we have the
following
Artin map
\[
\Phi_{F_n}:\ I_0 (F_n)\to \Gal (K_n/F_n),
\]
which induces an isomorphism $(I_0 (F_n)/H_n)_{p}\to \Gal (K_n/F_n)$.
Here $I_0 (F_n)$ is the group of ideals of $F_n$ which are prime to
$\la_n$, and $(I_0 (F_n)/H_n)_{p}$ is the
$p$-component of $I_0 (F_n)/H_n$.

The following facts were proved in \cite{K-M}:
\[
1)\,\,  \Gal^+  (K_n/F_n)\cong \Gal^+  (K_n/L_n)\cong\V_n^+
\]
\[
2)\,\,  K_{n-1}\cap F_n =F_{n-1}\, ({\mbox{lemma\,4.4}}).
\]
Obviously the field extension $F_n/F_{n-1}$ induces a natural
homomorphism
\[
\Gal  (K_{n-1}/F_{n-1})\cong (I_0 (F_{n-1})/H_{n-1})_{p} \to
(I_0 (F_n)/H_n)_{p}\cong \Gal  (K_n/F_n)
\]
which we denote with some abuse of notations by $\alpha_n$.
Therefore it is sufficient to prove that the latter $\alpha_n$ is
injective. First we note that the natural map
$F:\ \Gal  (K_{n-1}/F_{n-1})\to \Gal  (K_{n-1}F_n/F_{n})$ is an isomorphism.
Let us prove that $K_{n-1}F_n\subset K_n$.
Consider the Artin map
$\Phi_{F_n}':\ I_0 (F_n)\to \Gal (K_{n-1}F_n/F_n)$
(of course $F$ is induced by the canonical embedding $I_0 (F_{n-1})\to
I_0 (F_n)$).
We have to show that the kernel of  $\Phi_{F_n}'$ contains $H_n$.

To see this note that
$F^{-1}(\Phi_{F_n}' (s))=\Phi_{F_{n-1}}(N_{F_n/F_{n-1}}(s))$
for any $s\in I_0 (F_n)$. If $s\in H_n$ then
without loss of generality $s=1+\la_n^{p^n}t,\,\, t\in\Z[\zeta_{n}] $, and
thus,
$N_{F_n/F_{n-1}}(s))=1+pt_1$ for some $t_1\in\Z[\zeta_{n-1}]$. Now it
is clear that $\Phi_{F_n}' (s)=id_{K_n}$ since
$\Phi_{F_{n-1}}(1+pt_1)=id_{K_{n-1}}$.
%($0$ is the identical automorphism of $K_{n-1}$).

It follows that the identical map $id:\ I_0(F_n)\to I_0 (F_n)$
induces the canonical Galois surjection
$\Gal (K_n/F_n)\to \Gal (K_{n-1}F_n/F_n)$ and
we have the following commutative diagram:
\[
\xymatrix@=40pt{
& \Gal(K_{n-1}/F_{n-1})  \ar[dl]_{\alpha_n}
\ar[d]^{F}\\
\Gal(K_n/F_n)  \ar[r]         &    \Gal(K_{n-1}F_n/F_n)
}
\]
If $\alpha_n (a)=id$ then $F(a)=id$ and $a=id$ because $F$ is an
isomorphism which proves the lemma.
\end{proof}

\def\proofname{\bf Proof of Theorem~\ref{thm:Vn-semireg}.}
\begin{proof}
Induction with respect to $n$. If $n=1$ the result is known from for example \cite{K-M}.
Suppose the result holds with the index equal to $n-1$.
%There are no elements in $D_{n}^*$ with order greater than $p^n$ and hence there are no
%elements in $\V_n^+$ with order greater than $p^n$.
Proposition~\ref{prop:Vsurj} tells us that we have a surjection
$\pi_n:\V_n^+\rightarrow \V_{n-1}^+$
and Proposition~\ref{prop:pi-kernal}
that $\ker \pi_n$ isomorphic to
$(\Z/p\Z)^{r_{n-1}}$. Suppose $1+(x_{n-1}-1)^{k}$ is non-trivial in $\V_{n-1}^+$. Since
\begin{equation}\label{eq:ZDcomm2}
\xymatrix@=40pt{ \Z[\z_{n-1}]^{*+}
  \ar[r] \ar[d]^{\tilde{N}_{n,1}}  & D_{n}^{*+}  \ar[d]\\
  \Z[\z_{n-2}]^{*+} \ar[r] & D_{n-1}^{*+} }
\end{equation}
is commutative, $1+(x_{n}-1)^{k}$ is non-trivial in $\V_n^+$.
Moreover, since $\alpha_n$ is injective,
\[
\al(1+(x_{n-1}-1)^{k})=1+(x_{n}^p-1)^{k}=(1+(x_{n}-1)^{k})^p
\]
is non-trivial in $\V_n^+$. Now let $1+(x_{n-1}-1)^{s_i}$ generate $\V_{n-1}^+$ and suppose
$\pi_n(a_i)=1+(x_{n-1}-1)^{s_i}$. Since $\pi_n(1+(x_{n}-1)^{s_i})=1+(x_{n-1}-1)^{s_i}$
we get $a_i=b_i(1+(x_{n}-1)^{s_i})$ for some $b_i\in\ker \pi_n$, which implies that
$b_i^p$ is trivial.
Suppose $1+(x_{n-1}-1)^{s_i}$ has exponent $p^k$ for some $1\leq k \leq n-1$. To prove the
theorem we need to prove that $a_i$ has exponent $p^{k+1}$. Since
$\ker \pi_n\cong (\Z/p\Z)^{r_{n-1}}$, $a_i$ has exponent less than or equal to $p^{k+1}$.
But $(1+(x_{n-1}-1)^{s_i})^{p^k}=1+(x_{n-1}-1)^{p^k s_i}$ is non-trivial in $\V_{n-1}^+$
so
\[
a_i^{p^{k+1}}=b_i^{p^{k+1}} (1+(x_{n}-1)^{s_i})^{p^{k+1}}=(1+(x_{n}-1)^{s_i})^{p^{k+1}}
\]
is non-trivial in $\V_{n}^+$ by above, which is what we needed to show
\end{proof}
\def\proofname{\bf Proof.}
\section{A Weak Version of the Kervaire-Murthy Conjecture}\label{chap:weak}
In this section we will prove that $\cl \Q(\z_{n-1})(p)\cong \V_n^+/(\V_n^+)^p$.
Here $A(p):=\{x\in A:x^p=1\}$. It follows from Theorem~\ref{thm:Vn-semireg} that
$\V_n^+/(\V_n^+)^p$ has $r_{n-1}$ generators, and it was proved in [K-M] that
$ \cha \V_n^+$ can be embedded into  $\cl^{(p)} \Q(\z_{n-1})$.

So, in order to prove the result we need, it suffices to prove the following
\begin{thm}\label{thm:weakembedding}
There exists an embedding $\cl \Q(\z_{n-1})(p)\to \cha \V_n^+$.
\end{thm}
\begin{proof}
%XXX
First note that all our maps, $g_n, j_n, N_n$ etc and rings $A_n$ and can be extended $p$-adically.
Let $A_{n,(p)}$ be defined by
\[
A_{n,(p)}:=\frac{\Z_{p}[x]}{\Big(\frac{x^{p^{n}}-1}{x-1}\Big)},
\]
where $\Z_p$ denotes the ring of $p$-adic integers.
We have a commutative diagram
\begin{equation}\label{p-adicpullbackdiagram}
\xymatrix@=40pt{
A_{n,(p)} \ar[r]^{i_{n}} \ar[d]_{j_{n}}  & \Z_p[\z_{n-1}]
\ar[dl]^{N_{n-1}}
\ar[d]^{f_{n-1}}\\
A_{n-1,(p)} \ar[r]^{g_{n-1}}         &  D_{n-1}
}
\end{equation}
Let $U_{n,k,(p)}:=\{real\, \e\in \Z_p[\z_n]^*\, : \,\e\equiv 1 \m \la_n^k\}$
Considering pairs $(a, N_{n-1}(a))$, where $a\in \Z_p[\z_{n-1}]$,
we can embed $\Z_p[\z_{n-1}]^*$ into $A_{n,(p)}^*$. In [S2] it was proved
that $D_{n}^{*}$ is isomorphic to
$\Z_p[\z_{n-1}]^{*}/U_{n-1,p^{n}-1,(p)}$ (see also Lemma 2.6).
We hence have the following proposition
\begin{prop}\label{prop:Vnpadic}
%$\V_n^+$ is isomorphic to
\[
\V_n\cong\frac{\Z_p[\z_{n-1}]^{*}}{U_{n-1,p^{n}-1,(p)}\cdot g_n(\Z[\z_{n-1}]^*)}.
\]
\end{prop}
%YYY
Now for any valuation $\omega$ of $F_{n-1}=\Q (\z_{n-1})$ and any $a,b\in \Q (\z_{n-1})^*$
we have the norm residue symbol $(a,b)_{\omega}$ with values in the group
of $p$-th (not $p^n$) roots of unity. Let
$\omega=\la_{n-1}=(\z_{n-1}-\z_{n-1}^{-1})$
and
let $\eta_k=1-\la_{n-1}^k$. Then
\[
(\eta_i ,\eta_j )_{\la_{n-1}}= (\eta_i ,\eta_{i+j})_{\la_{n-1}}
(\eta_{i+j} ,\eta_{j})_{\la_{n-1}}(\eta_{i+j} ,\la_{n-1})_{\la_{n-1}}^{-j}
\]
It follows that $(a,b)_{\la_{n-1}}=1$ if $a\in U_{n-1,k},\,b\in U_{n-1,s}$ and
$k+s>p^n$. Further, $(\eta_{p^n} ,\la_{n-1})_{\la_{n-1}}=\z_0$ and therefore
$(\eta_i ,\eta_j )_{\la_{n-1}}\not= 1$ if $i+j=p^n$, $j$ is co-prime to $p$.

Let $\alpha$ be an ideal in $\Z[\z_{n-1}]$ co-prime to $\la_{n-1}$
and such that $\alpha^p=(q)$, where $q=1+\la_{n-1}^2 t\in \Z[\z_{n-1}]$
(we can choose such $q$
since $\z_{n-1}=1+\la_{n-1}\z_{n-1}(1+\z_{n-1})^{-1}$ and
 $\z_{n-1}(1+\z_{n-1})^{-1}\in\Z[\z_{n-1}]^*$).
Define the following action of $\cl \Q(\z_{n-1})(p)$
on $U_{n-1,2,(p)}^+$ :
\[
\tau_{\alpha}(v)=(v,q)_{\la_{n-1}}
\]
Let us prove that this action is well-defined. First of all it is
independent of the choice of the representative $\alpha$ in
$\cl \Q(\z_{n-1})(p)$ because if we use $r\alpha$ instead of $\alpha$
then $(v,r^p q)_{\la_{n-1}}=(v,q)_{\la_{n-1}}$.

The action is independent of the choice of $q$ by the following
reason: another generator of $\alpha^p$, which is $1$
modulo $\la_{n-1}^2$, differs from ``the old'' $q$ by some
unit $\gamma=1+\la_{n-1}^2 t_1$, and it can be easily verified that
$\gamma$ is either real or $\gamma=\z_{n-1}^{pk}\gamma_1 $ with
a real unit
$\gamma_1$. Hence we must  consider $\tau_{\gamma q}(v)$
for real $\gamma$. In other words we have to prove that
$(v,\gamma)_{\la_{n-1}}=1$. But if the latter is untrue, then
$(v,\gamma)_{\la_{n-1}}=\z_0$, which is not consistent with the
action of the ``complex conjugation'' ($v$ and $\gamma$ are real, while
$\z_0$ is not real).

Clearly $(U_{n-1,p^n-1,(p)},q)_{\la_{n-1}}=1$.
It remains to prove that $(\gamma , q)_{\la_{n-1}}=1$ for any
unit $\gamma$ and we
will obtain an action of $\cl \Q(\z_{n-1})(p)$ on $\V_n^+$.
For this consider a field extension $F_{n-1}(q^{1/p})/F_{n-1}$.
Since $(q)=\alpha^p$, it can remify in the $\la_{n-1}$ only.
Then clearly $(\gamma , q)_{\omega}=1$ for any $\omega\not=\la_{n-1}$
and it follows from the product formula that $(\gamma , q)_{\la_{n-1}}=1$.

Therefore $\cl \Q(\z_{n-1})(p)$ acts on $\V_n^+$ and obviously
$\tau_{\alpha\beta}=\tau_{\alpha}\tau_{\beta}$.

The last stage is to prove that any $\alpha\in\cl \Q(\z_{n-1})(p)$ acts
non-trivially on $\V_n^+$.
Let $(q)=\alpha^p$ and let $q=1+\la_{n-1}^k t$ with some $k>1$ and
$t$, co-prime to $\la_{n-1}$.

Let us prove that $k<p^n -1$.
Assume that $k>p^n -1$. Then the field extension $F_{n-1}(q^{1/p})/F_{n-1}$
is unramified. It is well-known that if $p$ is semi-regular, then
$F_{n-1}(q^{1/p})=F_{n-1}(\gamma^{1/p})$ for some unit $\gamma$.
Kummer's theory says that $q=\gamma r^p$ and then obviously
$\alpha=(r)$, i.e. $\alpha$ is a principal ideal. So, it remains to
prove that the case $k=p^n -1$ is impossible. For this consider
$\z_{n-1}$ and take into account that
$\z_{n-1}=1+\la_{n-1}\z_{n-1}(1+\z_{n-1})^{-1}$.
Then clearly it follows from the properties
of the local norm residue symbol  $(\, , \,)_{\la_{n-1}}$ that
$(\z_{n-1},q)_{\la_{n-1}}\not= 1$. On the other hand
$(\z_{n-1},q)_{\omega}=1$ for any $\omega\not= \la_{n-1}$ because
$\zeta_{n-1}$ is a unit and the extension $F_{n-1}(q^{1/p})/F_{n-1}$
is unramified in $\omega$. Therefore $(\z_{n-1},q)_{\la_{n-1}}= 1$ by
the product formula and the case $k=p^n -1$ is impossible and
$k<p^n -1$.

Now let us consider the cyclic subgroup of $\cl \Q(\z_{n-1})(p)$
generated by $\alpha$ and all the $q_i$ which generate all
$\alpha^{ps}$ for non-trivial $\alpha^s$ (i.e. $s$ is co-prime to $p$).
Let us choose that $q\in U_{n-1,k,(p)}$,
which has the maximal value of $k$.

Then $gcd (k, p)=1$ (otherwise consider $q(1-\la_{n-1}^{k/p})^p$).
Next we prove that $k$ is odd. If untrue, consider the following element
from our set of $\{ q_i\}$, namely $q/\sigma (q)$, where $\sigma$ is
the complex conjugation. Easy computations show that if $k$ is even
for $q$, then $q/\sigma (q)\in U_{n-1,s,(p)}$ with $s>k$.
On the other hand $q/\sigma (q)$ is in our chosen set of  $\{ q_i\}$
because it generates some ideal from the class of $\alpha^2$ since
$\cl \Q(\z_{n-1})(p)=\cl \Q(\z_{n-1})(p)^-$.
%$a=(1-\la_{n-1}^{p^n-k})(1+\la_{n-1}^{p^n-k})^{-1}$
%Easy computations show that
%$a=\eta_{p^n-k}^2 a_1$, where $a\in U_{n-1, p^n-k+1, (p)}$.
%From the above properties of the $(\, , \,)_{\la_{n-1}}$ it follows that
%$(a,q)_{\la_{n-1}}$ is a faithful $p$-root of unity.
%On the other hand for the ``complex
%conjugation'' $\sigma$ we have: $\sigma (a)=a^{-1}$,
%$\sigma (q)=q^{-1}r^p$ because $\cl \Q(\z_{n-1})(p)=\cl^- \Q(\z_{n-1})(p)$.
%Thus $(a,q)=(\sigma (a), \sigma(q))=\sigma (a,q)$, which is a contradiction.
Therefore we have proved that $k$ is odd. Then
$(\eta_{p^n-k}, q)\not=1$ and this means that   $\eta_{p^n-k}$ is a non-trivial
element of $\V_n^+$ for which $\tau_{\alpha}(\eta_{p^n-k})\not= 1$.

The theorem is proved.

%Assume that for some $(q)=\alpha^p$ (with
%$q=1+\la_{n-1}^2 t$) $(\eta_{2k},q)_{\la_{n-1}}=1$ for all $k$.

%$U_{n-1,2,(p)}^+/(U_{n-1,2,(p)}^+ )^p$ is generated by
%$\eta_{2k}$ with $ gcd (k,p)=1$, and $(\eta_{2k},\eta_{2m})_{\la_{n-1}}=1$.

\end{proof}
One of the Kervaire-Murthy conjectures was that $\cha\V_n^+\cong \cl^{(p)}\Q(\z_{n-1})$.
Now we partially solve this conjecture.
\begin{cor}
$\cl \Q(\z_{n-1})(p)\cong \V_n^+/(\V_n^+)^p\cong (\Z/p\Z)^{r_{n-1}}$ (see Section 2 for the
definition of $r_{n-1}$).
\end{cor}
\begin{proof}
It remains to prove the second isomorphism only,
which follows from Theorem 2.14.
\end{proof}
Now it is clear that the Assumption 2 from [H-S], which we used
there to describe
$\V_n^+$, is valid for any semi-regular prime.
\begin{cor}
Any unramified extension of $\Q (\z_{n-1})=F_{n-1}$ of degree $p$ is of the form
$F_{n-1}(\epsilon^{1/p})/F_{n-1}$, where $\epsilon$ is a unit
satisfying $\epsilon\equiv 1\m \la_{n-1}^{p^n +1}$.
\end{cor}
Now let us consider the Iwasawa module $T_p (\Q)$ as a $\Z_p$-module.
It is known from the Iwasawa theory that $T_p (\Q)\cong \Z_p^{\la}$
for semi-regular $p$, where $\la$ is the Iwasawa invariant for $p$
(see [W, Corollary 13.29]) and consequently $\cl^{(p)} (F_N)$
has $\la$ generators as an abelian group for big $N$
Therefore we obtain the following
\begin{cor}
There exists an integer $N$ such that
$r_k=\lambda$ for $k>N$.
Moreover, any unramified extension of $\Q (\z_{k})=F_{k}$ of degree $p$ is of the form
$F_{k}(\epsilon^{1/p})/F_{k}$, where $\epsilon\in\Z[\z_N ]^*$ is a unit
satisfying $\epsilon\equiv 1\m\la_{N}^{p^{N+1} +1}$.
\end{cor}

Finally we obtain Kummer's Lemma for semi-regular primes
\begin{cor}
Let a unit $\epsilon\in\Z[\z_{n-1}]^*$ satisfy
$\epsilon\equiv r^p \, {\mbox{mod}}\, \la_{n-1}^{p^n -1}$. Then
$\epsilon=\gamma^p\gamma_1$ with units $\gamma, \gamma_1$ and
$\gamma_1 \equiv 1\m \la_{n-1}^{p^n +1}$.
\end{cor}
\begin{proof}
If $\epsilon\equiv r^p \, {\mbox{mod}}\, \la_{n-1}^{p^n -1}$ then
$r^{-p}\epsilon\equiv 1\, {\mbox{mod}}\, \la_{n-1}^{p^n -1}$ and it
follows from the proof of the theorem that in fact
$r^{-p}\epsilon\equiv 1 \, {\mbox{mod}}\, \la_{n-1}^{p^n}$. Then the extension
$F_{n-1}({\epsilon}^{1/p})/F_{n-1}$ is unramified and therefore
by Corollary 3.4 $\epsilon=\gamma^p\gamma_1$,
where $\gamma_1\equiv 1 \, {\mbox{mod}}\, \la_{n-1}^{p^n +1}$. Clearly, then $\gamma$ is a unit.
\end{proof}

\section{The Kervaire-Murthy Conjecture}\label{chap:KM}

This section is devoted to the proof the following theorem.
\begin{thm}\label{thm:fullKM}
Let $p$ be a semi-regular prime. Then $\cha \V_n^+\cong \cl^{(p)} \Q(\z_{n-1})$.
\end{thm}

We start by defening the ray group $H_n'$ of $F_n$ by
\[
H_n'=\{(a)\subseteq F_n : a\equiv 1 \m \la_n^{p^n-1}\}.
\]
Let  $I_0 (F_n)$ be the group of all ideals of
$F_n$ prime to $\la_n$ and let $P_{0,n}$ be the
group of all principal fractional ideal of $F_n$ prime to $\la_n$.
Let $K_n'/F_n$ be the $p$-part of the ray extension associated to
$H_n'$. Then the Artin map
gives us an isomorphism
\[
\Gal(K_n'/F_n)\cong (I_n/H_n')_p.
\]
\begin{lemma}\label{lemma:Vnideal}
Let $p$ be a semi-regular prime.
Then $\V_n\cong P_{0,n-1}/H_{n-1}'$ and
$(P_{0,n-1}/H_{n-1}')^+=(I_{0} (F_{n-1})/H_{n-1}')_p^+$.
\end{lemma}
Before the proof, again recall that we define $\V_n$ by
\[
\V_n:=\frac{D_n^*}{\g_n(\Z[\z_{n-1}]^*)},
\]
where $\Z[\z_{n-1}]^*$ is embedded in $A_n$ using the map $a\mapsto (a,N_n(a))$. Our definition is equivalent
to the one in \cite{K-M} by Proposition~\ref{Vn=Vn}.
\begin{proof}
First recall that by Lemma~\ref{lemma:gker} $(a,N_n(a))\equiv 1 \m p$ in $A_n$ if and only if
$a\equiv 1 \m \la_{n-1}^{p^n-1}$ in $\Z[\z_{n-1}]$, that is, if and only if $(a)\in H_{n-1}'$. By just
counting elements it follows that for any $b\in A_n$, $b\equiv 1 \m (x-1)$ there exists $a\in \Z[\z_{n-1}]$
such that $b\equiv (a,N_n(a)) \m p$. This shows that $g_n:\Z[\z_{n-1}]\rightarrow \V_n$ induces a well defined,
bijective homomorphism
\[
G:P_{0,n}/H_n'\rightarrow \V_n.
\]
%since $g_n(a,N_n(a))$ trivial is equivalent to $g_n(a,N_n(a))\equiv g_n(\e,N_n(\e))$ which is equivalent to
%$g_n(a\e^{-1},N_n(a\e^{-1}))=1$ which is equivalent to $(a)\in H_{n-1}'$.
%
For the second equality note that if $\alpha$ represents an ideal in
$(I_0 (F_{n-1})/H_{n-1}')_p^+$ then
$\alpha^{\tau}=\alpha\dot (a)$ for some $(a)\in H_{n-1}'$. Since  $\alpha^{P^N}\in H_{n-1}'$ for some
$N$ we also get that  $\alpha^{P^N}$ is principal. Hence $\alpha$ represents an element of $(\cl^{(p)}(F_{n-1}))^+$
which is trivial by the semi-regularity condition which means $\alpha$ is
principal.
\end{proof}

\begin{cor}
$\V_n^+\cong \Gal^+(K_n'/F_n)$.
\end{cor}
We will also use the following lemma.
\begin{lemma}
$K_{n-1}'F_n\subseteq K_n \subset K_n'$.
\end{lemma}
\begin{proof}
The second inclusion is trivial. For the first, consider the commutative diagram
\[
\xymatrix@=40pt{
& \Z[\z_{n}]^* \ar[dl]_{N_{n}}
\ar[d]^{f_{n}}\\
A_{n}^* \ar[r]^{g_{n}}         &  D_{n}^*
}
\]
Recall that $\Z[\z_{n-1}]^*$ is mapped into $A_n$ and that $N_n$ acts as the usual norm $\tilde{N}_{n,1}$.
It follows that if $b\in H_n$, that is $b\equiv 1 \m \la_n^{p^n}$, then
$\tilde{N}_{n,1}(b)\equiv 1 \m \la_{n-1}^{p^n-1}$, that is $\tilde{N}_{n,1}(b)\in H_{n-1}'$.
To show the inclusion we have to show that
$\Phi_{F_n}(b)$ acts trivially on $K_{n-1}'F_n$ if $b\in H_n$.
But since the restriction of the Artin map
$\Phi_{F_n}$ to $K_{n-1}'$ is
$\Phi_{F_{n-1}}\circ \tilde{N}_{n,1}$ this follows from above.
\end{proof}
Now we want to extend some results of \cite{K-M}. Let $F=\bigcup F_n$ and
let $K$ be the maximal abelian $p$-extension of $F$ such that only
the prime $\la$ (the unique extension of $\la_n$ to $F$) ramifies in $K$.
Clearly $K_n$ and $K_n'$ are subfields of $K$.

Let $E=\bigcup \Z[\z_n]^*$ and $M_n=F(E^{1/p^n})$.
It was shown in \cite{K-M} that $M_n\subset K$. Set $M=\bigcup M_n$.
The group $\Gal (K/M)$ was described by Iwasawa, namely
$\cha \Gal(K/M)\cong S$, where $S$ is the direct limit of
$S_n:=\cl^{(p)} \Q (\z_n)$ with respect to the canonical embeddings
$\cl^{(p)} \Q (\z_n)\to \cl^{(p)} \Q (\z_{n+1})$.

Now, since $p$ is odd, $\Gal(K/F)=\Gal^+(K/F)\oplus \Gal^-(K/F)$ and it was
explained in \cite{K-M} that for the semi-regular primes
\[
\Gal(K/M)=\Gal^+(K/M)=\Gal^+(K/F)=\cha (S)
\]
It follows that $\Gal (M/F)=\Gal^-(K/F)$ and if we define $K^+$ to
be the subfield of $K$ fixed by $\Gal^-(K/F)$ then we see that
$K=K^+ M$ and $K^+\cap M=F$. Moreover,
\[
\Gal(K^+/F)=\Gal(K/M)=\Gal^+(K/M)=\Gal^+(K/F)=\cha (S):=G^+
\]

%Further let us define $KS_n$ subfields of $K^+$:
We also need the subfields $KS_n$ of $K^+$ defined as follows:
$S_n$ is a subgroup of $S$ and the latter group is dual to $G^+$.
Let $S_n^{\perp}\subset G^+$ be the subgroup annihilating
$S_n$ and let $KS_n\subset K^+$ be fixed by $S_n^{\perp}$.
Then obviously we have $\Gal(KS_n/F)=\cha (S_n)$, and since
$S$ is the direct limit of $S_n$ it follows that
$K^+ =\bigcup KS_n$.

Similarly, starting from the extensions $F_n\subset K_n$
and $F_n\subset K_n'$ we can define extensions of $F_n$,
namely $K_n^+$ and ${K_n'}^+$ such that
$\Gal (K_n^+/F_n)=\V_n^+$ and $\Gal ({K_n'}^+/F_n)=\V_{n+1}^+$.
%Clearly, since $K_n'F_{n+1}\sebseteq K_{n+1}$ it follows that
%${K_n'}^+ F_{n+1}=K_{n+1}^+$
%\begin{lemma}
%${K_n'}^+ F_{n+1}=K_{n+1}^+$
%\end{lemma}
%\begin{proof}
%Clearly, ${K_n'}^+ F_{n+1}$ is a class field over $F_{n+1}$
%associated with the ray group
%\end{proof}
Since $K_n\cap F=F_n$ (see \cite{K-M}) we have
$\Gal (K_n F/F)\cong \Gal (K_n/F_n)$ and consequently
\[
\Gal (K_n^+M/M)\cong\Gal (K_n^+ F/F)\cong \Gal (K_n^+/F_n)=\V_n^+
\]
\begin{lemma}
$K_n^+ F\subset KS_{n-1}$
\end{lemma}
\begin{proof}
It was proved in \cite{K-M} that the canonical surjection
\[
\Gal (K^+/F)=\Gal^+ (K/F)=\cha (S)\to\V_n^+
\]
factors through $\cha (S_{n-1})$ and hence $KS_{n-1}$
contains $K_n^+ F$.
\end{proof}
\begin{thm}
$K^+=\bigcup K_n^+F$
\end{thm}
\begin{proof}
It suffices to prove that $KS_n$ is contained in
$K_N^+ F$ for some big $N$.
Results of Section~\ref{chap:weak} imply that for big $N$
both groups $\cha (S_N)$ and $\V_N^+$ have $\la$ generators,
where $\la$ is the Iwasawa invariant. We also know that
$\cha (S_N)$ has $p^{\la N +\nu}$ elements and $\V_N^+$ has
$p^{\la N +\nu_1}$ elements, where $\nu,\, \nu_1$ do not
depend on $N$. Moreover, it follows from the structure
of $\V_N^+$ that any cyclic component of $\V_N^+$ has
$p^{N+\nu_i}$ elements where $\nu_i$ also do not depend
on $N$. Therefore every cyclic component of $\cha (S_N)$
has more than $p^{N+\nu_i}$ elements.

Now we want to compare the kernels of two canonical surjections,
$\cha (S_N)\to \V_{N+1}^+$ and $\cha (S_N)\to \cha (S_n)$.
The first kernel has $p^{\nu -\nu_1 -p}$ elements.
Each cyclic component of the second kernel has
$p^{N+\nu_i-n_i}$ elements where $S_n\cong \bigoplus \Z/p^{n_i}\Z$.
Therefore for big $N$ the first kernel is contained in the second
and we deduce that $KS_n\subset K_{N+1}^+F$.
\end{proof}

Let us construct a homomorphism $r:\,S_n\to \cha \V_{n+1}^+$.
Choose an element $b\in S_n$ and its representative $\be$, an
ideal in $\Z[\z_n]$ co-prime to $p$. Then
$\be^{p^k}=(q)$, where $q=1+\la_n^2 t\in \Z[\z_n]$
(see the proof of Proposition~\ref{prop:Vnpadic}).
Then we know from \cite{K-M} that $F(q^{1/p^k})\subset K$ and
more exactly $F(q^{1/p^k})\subset K_N^+M$ for some $N$.
Without loss of generality we can assume that $N\geq n+1$.
For any $v\in \V_N^+=\Gal (K_N^+M/M)$ define
$\tau_b (v)=v(q^{1/p^k})\cdot q^{-1/p^k}$ (Kummer's pairing).
Since by Lemma~\ref{lemma:Vn-inj} $\V_{n+1}^+$ is a subgroup of $\V_N^+$,
we can finally define $r(b)\in \cha \V_{n+1}^+$ as the
corresponding restriction of $\tau_b$.
\begin{thm}
$r$ is injective.
\end{thm}
\begin{proof}
Before we start proving the theorem we need the following
\begin{lemma}
Let $\al $ be an ideal of $\Z [\z_n]$ such that $\al^p=(q)$
with $q\equiv 1\m \la_n^2$. Then $F_n (q^{1/p})\subset K_n'$.
\end{lemma}
\begin{proof}
We have to prove that
$\Phi_{F_n}(a)=id$ on $F_n (q^{1/p})$
for any $a\equiv 1 \m \la_n^{p^{n+1}-1}$.
By the Reciprocity Law (see for instance \cite{C-F}, ch. VII)
it is true if $a$ is a local norm from the $\la_n$-adic
completion of $F_n (q^{1/p})$. The latter is equivalent to
that of $(a,q)_{\la_n}=1$
and this was established in the proof
of Theorem~\ref{thm:weakembedding}, where the corresponding local symbol was defined.
\end{proof}
We continue the proof of the theorem.
It is sufficient to prove that $r$ is injective on
the subgroup of elements of order $p$.
So, let $b^p=(q)$ where we can assume that $q\equiv 1\m \la_n^2$.
Then $q\in K_n'\subset K_{n+1}$. We have to find
\[
v\in \V_{n+1}^+ =\Gal^+ (K_{n+1}/F_{n+1})=
({{P_{0,n+1}}/{H_{n+1}}})^+
\]
such that
$v(q^{1/p})\times q^{-1/p}\neq 1$. So, without loss of generality
we can assume that $v\in \Z [\z_{n+1}]$
and $v\equiv 1 \m \la_{n+1}$. Further, $v$ acts on $K_{n+1}$
as $\Phi_{F_{n+1}} (v)$. On the other hand $q^{1/p}\in K_n'$
and therefore $\Phi_{F_{n+1}} (v) (q^{1/p})=
\Phi_{F_{n}} (N_{F_{n+1}/F_n}(v))(q^{1/p})$. By Proposition~\ref{Vn=Vn}
and Lemma~\ref{lemma:Vnideal}, $N_{F_{n+1}/F_n}$ induces an isomorphism between
${{P_{0,n+1}}/{H_{n+1}}}$ and $P_{0,n}/H_n'$. Thus we have to
find $w\in \Z[\z_n]$ such that $\Phi_{F_n}(w)(q^{1/p})\neq q^{1/p}$.
Again, by the Reciprocity Law we have that $\Phi_{F_n}(w)\psi_{\la_n}(w)=id$,
where $\psi_{\la_n}$ is the local Artin map. We get that
\[
\Phi_{F_n}(w)(q^{1/p})\times q^{-1/p}=\psi_{\la_n}(w^{-1})(q^{1/p})
\cdot q^{-1/p}=(q,w)_{\la_n}
\]
(the symbol $(\, , \, )_{\la_n}$ was defined  in the proof of Theorem~\ref{thm:weakembedding})
Then the required $w$ exists by Theorem~\ref{thm:weakembedding}.
\end{proof}
\begin{cor} Theorem~\ref{thm:fullKM} holds, i.e.
$\cl^{(p)} (F_n)\cong \cha \V_{n+1}^+$
\end{cor}
\begin{cor}
For any semi-regular prime $p$
\[
1)\ \ \cl^{(p)} \Q(\z_0)\cong \pic^{(p)} \Z C_p\cong (\Z/p\Z)^{r(p)}
\]
\[
2)\ \ \pic^{(p)} \Z C_{p^2}\cong (\Z/p\Z)^{{\frac{p-3}{2}}+r_1 -r(p)}
\oplus (\Z/p^2\Z)^{2r(p)}
\]
\end{cor}
%\begin{cor}
%Let $p$ be semi-regular and let $\e\in \Z[\z_0]^*$
%satisfy $\e\equiv 1 mod (p^2)$. Then $\e=\g^p$ for some
%$\g\in\Z[\z_0]^*$.
%\end{cor}

\section{Applications}\label{chap:appl}
{\bf The case $\la=r(p)$}
%\section{The Kervaire-Murthy conjectures when $r_n=r(p)$}

We now proceed by making an assumption under which we will prove
all of the  Kervaire-Murthy conjectures.
\begin{ass}\label{ass:r(p)r(p)r(p)}
$\la=r(p)$
%$\rank_p(\cl^{p}(\Q(\z_{n}))^-)=r(p)$ for all $n$.
\end{ass}
Then it follows that $\cl^{(p)} \Q(\z_{n}) \cong (\Z/p^{n+1}\Z)^{r(p)}$
for all $n$
%This holds for example
%if the Iwasava invariant $\la$
%and $\nu$
% satisfy $\la=r(p)=:r$
%in which case $\cl^{(p)}(\Q(\z_{n}))^-\cong (\Z/p^{n+1}\Z)^r$.
The assumption $\la=r(p)$ follows from
certain
congruence assumptions on
Bernoulli numbers (see page 202 in \cite{W})
known to hold for all primes
less than 4.000.000. Of course all these primes are semi-regular.

The following result
follows directly from Theorem~\ref{thm:Vn-semireg}.

\begin{thm}\label{thm:main2}
If $p$ is a semi-regular prime and $r$ the index of
irregularity and Assumption~\ref{ass:r(p)r(p)r(p)} holds, then
$\V_n^+ \cong (\Z/p^{n}\Z)^r$.
\end{thm}

We now proceed to show how we can directly show that $\V_n^+ =V_n^+$ when $\V_n^+ \cong (\Z/p^n \Z)^r$.
The proof of
this relies of constructing a certain basis for $D_{n-1}^{+}$
consisting of norms of elements from $\Z[\z_{n-1}]^*$
considered$\mod p$.

Let $\Phi:U_{n-1,p^n-p^{n-1}}\rightarrow D_{n-1}^+$ be defined
  by
\[
\Phi(\e)=N_{n-1}\big( \frac{\e -1}{p}\big)-\frac{N_{n-1}(\e)-1}{p}
\mod p.
\]
Since $N_{n-1}$ is additive$\mod p$ one can show with some simple calculations
that $\Phi$ is a group homomorphism. See Lemmas~\ref{lemma:varphihomo} and~\ref{lemma:omegahomo} for details.

Explicitly, what we want to prove is the following.
\begin{thm}\label{thm:Phihomo}
If $\V_n^+ \cong (\Z/p^n \Z)^r$, then $\Phi$ is a surjective group homomorphism.
\end{thm}
%\begin{thm}
%If $\V_n^+ \cong (\Z/p^n \Z)^r$ then
%there are elements $\e_k \in U_{n-1,p^n-p^{n-1}}\rightarrow D_{n-1}^+$, $k=1,\hdots, p^{(p^{n-1}-1)/2}$ such that
%$\Phi(\e_k)$
%\[
%N_{n-1}\big( \frac{\e_k -1}{p}\big)-\frac{N_{n-1}(\e_k)-1}{p}
%\mod p
%\]
%forms a basis for $D_{n-1}^{+}$ over $\F_p$.
%\end{thm}
As we can see by the following corollary, the theorem is what we need.
\begin{cor}\label{cor:v=v}
If $\V_n^+ \cong (\Z/p^n \Z)^r$, then $V_n^+ = \V_n^+$
\end{cor}
\def\proofname{\bf Proof of the Corollary.}
\begin{proof}
We want to show
that
for any $(1,\ga)\in A_{n}^*$ there exists $(\e,N_{n-1}(\e)) \in A_{n}^*$ such that
$(1,\ga) \equiv (\e,N_{n-1}(\e)) \mod p$, or more explicitly
that for
all $\ga \in A_{n-1}^{*+}$, $\ga\equiv 1 \m p$ there exists $\e \in \Z[\z_{n-1}]^*$ such
that $(\e,N(\e))\equiv (1,\ga) \mod p$ in $A_{n}$. This is really
equivalent to the following three statements in $\Z[\z_{n-1}]$,
$A_{n-1}$ and $D_{n-1}$ respectively
\begin{eqnarray*}\label{eqn:3congs}
\e                          &\equiv& 1 \mod p\\
N_{n-1}(\e)                       &\equiv& \ga \mod p\\
N_{n-1}\big( \frac{\e -1}{p}\big) &\equiv& \frac{N_{n-1}(\e)-\ga}{p} \mod p
\end{eqnarray*}
%$\e \equiv 1 \mod p$ implies that $N_{n-1}(\e)\equiv 1 \mod p^2$, so
%$\frac{N_{n-1}(\e)-\ga}{p} \equiv \frac{1-\ga}{p} \mod p$.
Note that $(1,\ga)\in A_{n}$ implies $g_{n-1}(\ga)=f_{n-1}(1)$
in $D_{n-1}$, or in other words, that $\ga\equiv 1 \mod p$.
Hence we only need to show that for any $\ga \in A_{n-1}^{*+}$ there exists
$\e \in U_{n-1,p^n-p^{n-1}}$ such that
\[
N_{n-1}\big( \frac{\e -1}{p}\big) - \frac{N_{n-1}(\e)-1}{p}
\equiv \frac{1-\ga}{p} \mod p.
\]
But the left hand side is exactly $\Phi(\e)$ so the corollary really does follow from Theorem~\ref{thm:Phihomo}
\end{proof}
\def\proofname{\bf Proof.}
%%%%%%%%%%%%%%%%%%%%%%%%%%%%%start proof thm:Phihomo
We now proceed to start proving Theorem~\ref{thm:Phihomo}.
Recall that $r=r(p)$ are the number of indexes $i_1,i_2\hdots i_r$
among $1,2\hdots (p-3)/2$ such that the nominator of the Bernoulli
number $B_{i_k}$ (in reduced form) is divisible by $p$.

Let $\bar{E_n}:D_{n}\rightarrow D_{n}^*$ be the
truncated exponantial map defined by
\[
\bar{E_n}(y)=1+y+\frac{y^2}{2!}+\hdots +\frac{y^{p-1}}{(p-1)!}
\]
and let $\bar{L_n}:D_{n}^*\rightarrow D_{n}$ be the truncated logarithm map
\[
\bar{L_n}(1+y)=y-\frac{y^2}{2}+\hdots -\frac{y^{p-1}}{(p-1)}.
\]
We also consider the usual $\la$-adic $\log$-map defined by a power series as usual.
%For simplicity, let the ideal $la_0$ be represented by $(\z_0-\z_0^{-1})$.

%Recall from the proof of Lemma~\ref{lemma:Nsurj} that
We denote the cyclytomic units
of $\Z[\z_0]^{*+}$ by $C_0^+$.
Let $\mathcal{M}$ be the group of real $\la_0$-adic integers with zero trace. Any
$a\in \mathcal{M}$ can be uniquely presented as $a=\sum_{i=1}^{m-1} b_i \la_0^{2i}$,
$m=(p-1)/2$. Consider the homomorphism $\Psi:\Z[\z_0]^*\rightarrow \mathcal{M}$ defined by
$\e\mapsto \log(\e^{p-1})$. Following \cite{B-S}, page 370-375, we see that there are exactly
$r$ elements $\la_0^{2i}$, namely $\la_0^{2i_k}$, such that $\la_0^{2i}\not\in \Psi(C_0^+)$.
This implies
that for exactly the $r$ indexes $i_1,i_2\hdots i_r$ we have
$(\bar{x}_{1}-\bar{x}_{1}^{-1})^{2i_k}\neq g_{1}(log(\e^{p-1}))$ for any $\e\in C_0^+$.

Suppose
$(x-x^{-1})^{2i_s}=g_{1}(\log \e)$ for some $\e \in \Z[\z_0]^{*+}$.
%By the same argument as in the proof of Lemma~\ref{lemma:Nsurj}
It is well known that the index of $ C_0^+$ in $\Z[\z_0]^{*+}$ equals the classnumber
$h_+$ of $\Q(\z_0)^+$. Since $p$ is semi-regular there exists $s$ with $(s,p)=1$ such that
$\e^{s}\in C_0^+$ and by co-primality of $s(p-1)$ and $p$ we can find $u,v$ such that
$1=s(p-1)u+pv$. Then $\e=\e^{s(p-1)u+pv}=(\e^s)^{p-1}\e^{pv}$ so
$\log ((\e^{su})^{(p-1)u})=\log\e-pv\log\e \equiv \log\e \equiv (x-x^{-1})^{2i_s}$, which is a contradiction.
Hence $(x-x^{-1})^{2i_s}\not\in g_{1}(\log \Z[\z_0]^{*+})$. Since formally, $\exp(\log(1+y))=1+y$ it is not
hard to see that $E_1(L_1(1+y))\equiv 1+y\mod p$ and that we have a commutative diagram
\[
\xymatrix@=40pt{
 \tilde{\Z}[\z_{0}]^{*+}
  \ar[d]_{\log} \ar[dr]_{L_{1}}  \ar[drr]^{g_{1}} \\
  \mathcal{M} \ar[r]^{\mod p} & D_{1}^{+} \ar@<1ex>[r]^{\bar{E}_1} & D_{1}^{*+} \ar@<1ex>[l]^{\bar{L}_1}
}
\]
Recall that $D_{n,(s)}^{*+}:=\{y\in D_{n}^{*+} : y\equiv 1 \mod
(x-x^{-1})^s \}$ and that
we know that $\V_1^+:=D_{1}^{*+}/g_{1}(\Z[\z_0]^{*+})$ has $r:=r(p)$ generators.
If we now apply the map $E_1$ and do some simple calculations we now get the following proposition.
\begin{prop}\label{prop:V1-generators}
  The $r$ elements $\bar{E_1}((x_{1}-x_{1}^{-1})^{2i_k})$ generate
  $D_{1}^{*+}/g_{1}(\Z[\z_0]^{*+})$ and belong to
  $D_{1,(2)}^{*+}$ but do not belong to $D_{1,(p-2)}^{*+}$.
\end{prop}

We now want to lift this result to $D_{n}^{*+}$. From now on (exepting Lemma 4.11)
we will denote the generator $x\in D_n$ by $x_n$.
\begin{prop}\label{prop:Dn-Zn-generators}
  If Assumption 1 holds, then the $r$ elements
  $\bar{E_n}((x_{n}-x_{n}^{-1})^{2i_k})$ generate the group
  $\V_n^+:=D_{n}^{*+}/g_{n}(\Z[\z_{n-1}]^{*+})$. The elements
  $\bar{E_n}((x_{n}-x_{n}^{-1})^{2i_k})^{p^{n-1}}$
  are non-trivial in $\V_n^+$,
  belong to
  $D_{n,(p^{n-1})}^{*+}$ but do not belong to
  $D_{n,(p^n-2p^{n-1})}^{*+}$
\end{prop}
\begin{proof}
  Induction on $n$. If $n=1$ this is exactly
  Proposition~\ref{prop:V1-generators}.  Suppose the statement holds
  for the index equal to $n-1$.  The diagram
\begin{equation}\label{eq:ZDcomm}
\xymatrix@=40pt{ \Z[\z_{n-1}]^{*+}
  \ar[r] \ar[d]^{\tilde{N}_{n,1}}  & D_{n}^{*+}  \ar[d]\\
  \Z[\z_{n-2}]^{*+} \ar[r] & D_{n-1}^{*+} }
\end{equation}
is commutative.  Hence, if $z_n \in D_{n}^*$ is mapped to $z_{n-1}
\in D_{n-1}^*$ and $z_{n-1} \not\in \Image \Z[\z_{n-2}]^{*}$, then
$z_n \not\in \Image \Z[\z_{n-1}]^{*}$. Moreover, $z_n^p\not\in \Image
\Z[\z_{n-1}]^{*}$ in this case.
This follows from the fact that $\V_{m}^+\cong (\Z/p^{m}\Z)^r$ for all $m$.
Hence, if an element
$z \in \V_{n}^+$ has order $p$, then the surjection $\V_{n}^+\rightarrow \V_{n-1}^+$
maps $z$ to the neutral element in
$\V_{n-1}^+$.
%%%%%%%%%%%%%%%%%
%To prove this latter statement we will use that
%$\V_n^+\cong (\Z/p^n\Z)^r$, $\V_{n-1}^+\cong (\Z/p^{n-1}\Z)^r$ and the
%surjection $D_{n}\rightarrow D_{n-1}$ induces a surjection
%$(\Z/p^n\Z)^r\rightarrow (\Z/p^{n-1}\Z)^r$ If $a\neq 0$ in
%$(\Z/p^n\Z)$ but $a^p=0$, then $a$ maps to $0$ in $\Z/p^{n-1}\Z$. This
%translates to the fact that $z_n^p\not\in \Image \Z[\z_{n-1}]^{*}$,
%which is what we wanted to show.
%%%%%%%%%%%%%%%%%%%%%%%%%%%
Now, the elements
$\bar{E_{n}}((x_{n}-x_{n}^{-1})^{2 i_k})^{p^{n-1}}$
%=1+(x_{n}-x_{n}^{-1})^{2 p^{n-1}i_k}+\hdots
%\in D_{n}^{*+}(2p^{n-2})\setminus D_{n}^{*+}(p^{n-1}-2p^{n-1})$
are not in the image of $\Z[\z_{n-1}]^{*}$ by Theorem 4.3 and
since $\bar{E_{n}}((x_{n}-x_{n}^{-1})^{2 i_k})^{p^{n-2}}$
clearly map onto
$\bar{E}_{n-1}((x_{n-1}-x_{n-1}^{-1})^{2i_k})^{p^{n-2}} \not\in g_{n-1}(\Z[\z_{n-2}]^{*+})$
by induction.  Finally, since $1\leq
2i_k\leq p-1$ we get $p^{n-1} \leq 2p^{n-1}i_k\leq p^{n}-2p^{n-1}$ and
this means that all the elements
\begin{eqnarray*}
\bar{E_{n}}((x_{n}-x_{n}^{-1})^{2 i_k})^{p^{n-1}}&=&
(1+(x_{n}-x_{n}^{-1})^{2 i_k}+\hdots)^{p^{n-1}}=\\
&=&1+(x_{n}-x_{n}^{-1})^{2p^{n-1} i_k}+\hdots
\end{eqnarray*}
fulfil our requirements.
\end{proof}

%\subsection{The function $\varphi$}
Recall that $c:D_{n}\rightarrow D_{n}$ is the map induced by
$\bar{x} \mapsto \bar{x}^{-1}$ and that
$D_{n}^+:=\{a\in D_{n}:c(a)=a\}$
Define $\varphi:U_{n-1,p^n-p^{n-1}}^{+}\rightarrow D_{n-1}^+$ by
$\varphi(\ga)=N_{n-1}\big( \frac{\ga -1}{p}\big) \mod p$.
\begin{lemma}\label{lemma:varphihomo}
$\varphi$ is a homomorphism from the multiplicative group
$U_{n-1,p^n-p^{n-1}}^{+}$
to the additive group $D_{n-1}^+$ and the kernel is
$U_{n-1,p^n-1}^{+}=U_{n-1,p^n+1}^{+}$.
\end{lemma}
\begin{proof}
Let $\e$ and $\ga$ belong to $\in U_{n-1,p^n-p^{n-1}}^{+}$. Then, since $N_{n-1}$ is additive mod $p$
and $N_{n-1}(\e)\equiv 1 \mod p$,
\begin{eqnarray*}
N_{n-1}\big(\frac{\e\ga -1}{p}\big) & \equiv & N_{n-1}\big(\frac{\e(\ga -1)+ (\e -1)}{p}\big)\equiv \\
& \equiv & N_{n-1}(\e)N_{n-1}\big(\frac{\ga -1}{p}\big)+N_{n-1}\big(\frac{\e -1}{p}\big)\equiv \\
& \equiv & N_{n-1}\big(\frac{\ga -1}{p}\big)+N_{n-1}\big(\frac{\e -1}{p}\big) \mod p
\end{eqnarray*}
so $\varphi$ is a homomorphism. Suppose $N_{n-1}((\ga -1)/p)\equiv 0 \mod p$. Then, by
Proposition~\ref{prop:Norms},
$f_{n-1}((\ga -1)/p)=0$ which means $\ga \in U_{n-1,p^n-1}^{+}=U_{n-1,p^n+1}^{+}$ (the latter
equality is due to Lemma 3.2).
%By
%INSERT APROPRIATE REFERENCE %~\ref{ST2 p^n-1 => p^n}
%and
%INSERT APROPRIATE REFERENCE %~\ref{p^n => p^n+1}
%$U_{n-1,p^n-1}^{+}=U_{n-1,p^n+1}^{+}$.
\end{proof}

In this notation, what we want to prove is the following
\begin{prop}\label{prop:iso}
If Assumption 1 holds, then the map
\[
\tilde{\varphi}:(U_{n-1,p^n-p^{n-1}})/(U_{n-1,p^n+1})
\rightarrow D_{n-1}^+
\]
induced by $\varphi$ is an isomorphism.
\end{prop}
Since $\tilde{\varphi}$ is obviously injective it is enough to prove the
following proposition
\begin{prop}\label{prop:U/U-order}
Suppose Assumption 1 holds. Then
\[
|D_{n-1}^+|=\Big|\frac{U_{n-1,p^n-p^{n-1}}}{U_{n-1,p^n+1}}\Big|.
\]
\end{prop}
\begin{proof}
Recall that $|D_{n-1}^+|=p^{\frac{p^{n-1}-1}{2}}$ so we need to prove that
\[
|(U_{n-1,p^n-p^{n-1}})/(U_{n-1,p^n-1})|=p^{\frac{p^{n-1}-1}{2}}.
\]
An element of $\V^+_{n}$ of the form $b=1+(x_n-x_n^{-1})^{2s_1}$, where
$p^{n-1}<2s\leq 2s_1<p^n -1$,
correspond to a non-trivial element of
\[
\frac{D_{n,(2s)}^{*+}}{g_{n}(\Z[\z_{n-1}]^{*+})\cap D_{n,(2s)}^{*+}}
\]
which is a canonical subgroup of $\V_n^+$.
%since $U_{n-1,p^{n}-1}=\kernel ({g_{n}}{|_{\Z[\z_{n-1}]^*}})$.
If $t_{2s}$ is the number of independent such elements $b$, then
\[
\frac{D_{n,(2s)}^{*+}}{g_{n}(\Z[\z_{n-1}]^{*+})\cap D_{n,(2s)}^{*+}}\cong(\Z/p\Z)^{t_{2s}}
\]
By Proposition~\ref{prop:Dn-Zn-generators}, $t_{2s}=0$ if $2s>p^n-2p^{n-1}$.
On the other hand \[g_{n}(\Z[\z_{n-1}]^{*+})\cap D_{n,(2s)}^{*+}\cong U_{n-1,2s}/U_{n-1,p^n -1}\]
since $U_{n-1,p^{n}-1}=\kernel (g_{n})$,
and hence $U_{n-1,2s}/U_{n-1,p^n -1}\cong D_{n,(2s)}^{*+}$ if $2s>p^n-2p^{n-1}$.
The number of elements in $D_{n,(2s)}^{*+}$ is $p^{\frac{p^n -1 -2s}{2}}$.
Setting $2s=p^{n}-p^{n-1}$ completes the proof.
\end{proof}

%\subsection{Fine structure of the norm stuff}
%We start by proving the following lemma
We now have to do carefull estimations of some congruences of our norm-maps.

\begin{lemma}\label{lemma:normorder}
Let $2\leq n$ and $1 \leq k < n$. If $\e\in\Z[\z_{n-1}]$ and
If $\e\equiv 1 \mod p^{s+1}\la_{n-1}^{p^{n-1}-p^{k}}$, then
$(N_{n-1}(\e)-1)/p$ can be represented by a polynomial $f(x)=p^sf_1(x)$ in
$A_{n-1}$, where $f_1(x)\equiv 0 \mod (x-1)^{p^{n-1}-p^{k-1}}$ in  $D_{n-1}$.
\end{lemma}
Before the proof, recall that the usual norm $\tilde{N}_{n,1}$, $1 \leq n, 1\leq k < n$,
can be viewed as
a product of automorphisms of $\Q(\z_n)$ over $\Q(\z_{n-1})$. If
$t_n \in \Z[\z_n]$ and  $t_{n-1} \in \Z[\z_{n-1}]$ we immediately get
$\tilde{N}_{n,1}(1+t_{n-1}t_n)=1+\trace_{\Q(\z_n)/\Q(\z_{n-1})}(t_n)t_{n-1}t'$
for some $t'\in \Z[\z_{n-1}]$. Recall that the trace is always divisible by $p$.
In the proof below we will for convenience denote any generic element whose value
is not interesting for us by the letter $t$.
\begin{proof}
Induction on $n$. If $n=2$ (which implies $k=1$),
$N_{n-1}=\tilde{N}_{1,1}:\Z[\z_1]\rightarrow A_{1}\cong\Z[\z_0]$.
Let $\e:=1+tp^{s+1}$ Then
$\e=1+tp^{s}\la_1^{p^2-p}=1+tp^{s}\la_0^{p-1}$. By the note above,
\[
\frac{\tilde{N}_{1,1}(\e)-1}{p}=tp^s\la_0^{p-1}
\]
which is represented by some $f(x)=p^s(x-1)^{p-1}f_1(x)$ in $A_{1}$
Suppose the statement of the Lemma holds with the index equal to $n-2$. Let
$\e:=1+tp^{s+1}\la_{n-1}^{p^{n-1}-p^k}$. Note that $\e=1+tp^{s+1}\la_{n-2}^{p^{n-2}-p^{k-1}}$
and by the note before this proof, $\tilde{N}_{n-1,1}(\e)=1+tp^{s+2}\la_{n-2}^{p^{n-2}-p^{k-1}}$.
Let $(N_{n-1}(\e)-1)/p$ be represented by a pair $(a,b)\in \Z[\z_{n-2}]\times A_{n-2}$.
Then $a=(\tilde{N}_{n-1,1}(\e)-1)/p=tp^{s+1}\la_{n-2}^{p^{n-2}-p^{k-1}}$.
In $A_{n-2,1}$ $a$ hence can be represented by a polynomial
$a(x)=p^{s+1}(x-1)^{p^{n-2}-p^{k-1}}a_1(x)$ for some $a_1(x)$.
By the expression for
$\tilde{N}_{n-1,1}(\e)$ and by the assumption, we get
\[
b=\frac{N_{n-2}(\tilde{N}_{n-1,1}(\e))-1}{p}=
\frac{N_{n-2}(1+tp^{s+2}\la_{n-2}^{p^{n-2}-p^{k-1}})-1}{p}=p^{s+1}b_1(x)
\]
where $b_1(x)\equiv (x-1)^{p^{n-2}-p^{k-2}}b_2(x) \mod p$ for some $b_2(x)$.
Define $b(x):=p^{s+1}b_1(x)$.
We want to find a polynomial $f(x)\in A_{n-1}$ that represents $(a,b)$, that is,
maps to $a(x)$ and $b(x)$ in $A_{n-2,1}$ and $A_{n-2}$ respectively.
Note that
\[
p=\frac{x^{p^{n-1}}-1}{x^{p^{n-2}}-1}+t(x)\frac{x^{p^{n-2}}-1}{x-1}
\]
for some polynomial $t(x)\in \Z[x]$. Hence
\[
a(x)-b(x)=\Big(\frac{x^{p^{n-1}}-1}{x^{p^{n-2}}-1}+t(x)\frac{x^{p^{n-2}}-1}{x-1} \Big)
p^s((x-1)^{p^{n-2}-p^{k-1}}a_1(x)-b_1(x))
\]
Then we can define a polynomial $f(x)$ by
\begin{eqnarray*}
f(x):&=&a(x)+p^s((x-1)^{p^{n-2}-p^{k-1}}a_1(x)-b_1(x))\frac{x^{p^{n-1}}-1}{x^{p^{n-2}}-1}=\\
&=&b(x)+p^s((x-1)^{p^{n-2}-p^{k-1}}a_1(x)-b_1(x))t(x)\frac{x^{p^{n-2}}-1}{x-1}.
\end{eqnarray*}
Clearly, $f$ maps to $a(x)$ and $b(x)$ respectively.
We now finish the proof by observing that
\begin{eqnarray*}
f(x)/p^s&=&p(x-1)^{p^{n-2}-p^{k-1}}a_1(x)+((x-1)^{p^{n-2}-p^{k-1}}a_1(x)-b_1(x))\frac{x^{p^{n-1}}-1}{x^{p^{n-2}}-1}\equiv\\
&\equiv&((x-1)^{p^{n-2}-p^{k-1}}a_1(x)-(x-1)^{p^{n-2}-p^{k-2}}b_2(x))(x-1)^{p^{n-1}-p^{n-2}}=\\
&=&(a_1(x)-(x-1)^{p^{k-1}-p^{k-2}}b_2(x))(x-1)^{p^{n-1}-p^{k-1}} \mod p.
\end{eqnarray*}
\end{proof}

By setting $s=0$ we in the lemma above we immediately get the following theorem.
\begin{thm}\label{thm:varphiorder}
Let $2\leq n$ and $1 \leq k < n$.  Suppose
$\e \in U_{n-1,p^n-p^k}$.
%$\e \equiv 1 \mod \la_{n-1}^{p^{n}-p^{k}}$ in $\Z[\z_{n-1}]$.
Then $g_{n-1}((N_{n-1}(\e)-1)/p) \equiv 0 \mod (x-1)^{p^{n-1}-p^{k-1}}$
in $D_{n-1}$
\end{thm}

The following proposition is immediate by using that
$g_{n-1}N_{n-1}=f_{n-1}$.
\begin{prop}\label{prop:omegaorder}
Let $2\leq n$, $1 \leq k < n$  and let
$\e \in U_{n-1,p^n-p^k}\setminus U_{n-1,p^n-p^{k-1}}$.
%$\e \equiv 1 \mod \la_{n-1}^{p^{n}-p^{k}}$ in $\Z[\z_{n-1}]$.
Then $g_{n-1}((N_{n-1}((\e-1)/p)))
\equiv 0 \mod (x-1)^{p^{n-1}-p^{k}}$ but
$g_{n-1}((N_{n-1}((\e-1)/p)))
\not\equiv 0 \mod (x-1)^{p^{n-1}-p^{k-1}}$
in $D_{n-1}$.
\end{prop}
Let $\omega:U_{n-1,p^n-p^{n-1}}\rightarrow D_{n-1}^+$ be defined by
$\omega(\ga):=g_{n-1}((N_{n-1}(\ga)-1)/p)$.
\begin{lemma}\label{lemma:omegahomo}
$\omega$ is a homomorphism
\end{lemma}
\begin{proof}
Suppose $\e$ and $\ga$ belong to $U_{n-1,p^n-p^{n-1}}$. Then
$N_{n-1}(\ga)\equiv 1 \mod p$ in $A_{n-1}$ because
\[
N_{n-1}(\ga)=(\tilde{N}_{n-1,1}(\ga),\tilde{N}_{n-1,2}(\ga),\hdots,\tilde{N}_{n-1,n-1}(\ga))
\]
and $\tilde{N}_{n-1,k}(\ga)\equiv 1 \mod p^2$ for all $k=1,2,\hdots,n-1$.
Hence
\begin{eqnarray*}
\omega(\e\ga)&\equiv&\frac{N_{n-1}(\e\ga)-1}{p}=
\frac{N_{n-1}(\ga)N_{n-1}(\e)-N_{n-1}(\e)+N_{n-1}(\e)-1}{p}\equiv\\
&\equiv& N_{n-1}(\ga)\frac{N_{n-1}(\e)-1}{p}+\frac{N_{n-1}(\ga)-1}{p}\equiv\\
&\equiv& \frac{N_{n-1}(\e)-1}{p}+\frac{N_{n-1}(\ga)-1}{p}=\omega(\e)+\omega(\ga) \mod p
\end{eqnarray*}
\end{proof}

Note that if $\e\in U_{n-1,p^n-1}$ then $\omega(\e)=0$. This can be shown using similar,
but simpler,
methods as we did in the proof of Lemma~\ref{lemma:normorder}. We can hence define
\[
\tilde{\omega}:\frac{U_{n-1,p^n-p^{n-1}}}{U_{n-1,p^n-1}}\rightarrow D_{n-1}^+.
\]
Now, if $a\in D_{n-1}^+$, let ${\mathcal O}(a)$ be the maximal power of $(x-x^{-1})$
that devides $a$. In this language we can combine Thereom~\ref{thm:varphiorder} and
Proposition~\ref{prop:omegaorder} to the following lemma.
\begin{lemma}\label{lemma:omagavarphiorders}
Let $2\leq n$, $1 \leq k < n$  and let
$\e \in U_{n-1,p^n-p^k}\setminus U_{n-1,p^n-p^{k-1}}$. Then
$p^{n-1}-p^k \leq {\mathcal O}(\tilde{\varphi}(\e)) < p^{n-1}-p^{k-1} \leq
{\mathcal O}(\tilde{\omega}(\e))$.
\end{lemma}

\begin{prop}\label{prop:tildephiiso}
The map $\tilde{\Phi}:=\tilde{\varphi}-\tilde{\omega}$ is an isomorphism.
\end{prop}
\begin{proof}
By Proposition~\ref{prop:iso}
${\tilde{\varphi}}$ is an isomorphism. Hence there exists (classes of) units
$\e_i, i=1,2,\hdots,(p^{n-1}-1)/2$
such that the set ${\tilde{\varphi}}(\e_i)$ forms a basis for
$D_{n-1}^+$. If $a\in D_{n-1}^+$ there exist unique $a_i$ such that
$a=\sum_{i=1}^{(p^{n-1}-1)/2}a_i {\tilde{\varphi}}(\e_i)$.
To prove the Proposition it is enough to
show that the map
\[
\sum_{i=1}^{(p^{n-1}-1)/2}a_i {\tilde{\varphi}}(\e_i)\mapsto
\sum_{i=1}^{(p^{n-1}-1)/2}a_i ({\tilde{\varphi}}(\e_i)-{\tilde{\omega}}(\e_i))
\]
is invertible. Consider the matrix $M$ for this map in the basis
$\{(x-x^{-1})^{2j}\}$.
Obviously this matrix can be written $I-M'$, where $I$ is the
identity matrix and $M'$ is induced
by ${\tilde{\varphi}}(\e_i)\mapsto {\tilde{\omega}}(\e_i)$. By
Lemma~\ref{lemma:omagavarphiorders} the matrix $M'$
is a lower triangular matrix with zeros on the diagonal.
This means $M$ is lower triangular with ones on the
diagonal and hence invertible.
\end{proof}
\def\proofname{\bf Proof of Theorem~\ref{thm:Phihomo}.}
\begin{proof}
The map $\tilde{\Phi}$ is obviously induced by $\Phi$ which hence must be surjective by
prop~\ref{prop:tildephiiso}.
\end{proof}
\def\proofname{\bf Proof.}
%%%%%%%%%%%%%%%%%%%%%%%%%%%%%%%%%%%%%%%%%end proof thm:Phihomo
By Iwasawas theorem,
there are numbers $\la\geq 0$, $\mu \geq 0$ and $\nu$ such that
\[
|\cl^{(p)} \Q(\z_{n-1}^-|=p^{\la (n-1) + \mu p^n +\nu}
\]
for all $n$ big enough.
It has later been proved that $\mu =0$, so for big $n$,
$|\cl^{(p)} \Q(\z_{n-1})^-|=p^{\la (n-1) +\nu}$.
%It turns out that the Iwasawa invariant $\la$ is
%related to the sequence $\{r_n\}$ which we introduced in Chapter~\ref{chap:prel}.
%
%Recall that  \cite{K-M} provide us with an
%injection of $\cha(\V_n^+)$  into $\cl^{p}(\Q(\z_{n-1}))^-$,
% the $p$-component of the class group of $\Q(\z_{n-1})$. This shows that the number of elements in
%$\V_n^+$ is bounded by the number of elements in $\cl^{p}(\Q(\z_{n-1}))^-$.
%Since we know that $|\V_n^+|=p^ {r_o+r_1+\hdots +r_{n-1}}$
%this immediately implies the following proposition.
%\begin{prop}\label{prop:Vnorderbig}
%There is a number $n_0$ such that for $n\geq n_0$,
%$|\V_n^+|\leq p^{\la (n-1)+\nu}$
%\end{prop}
%By comparing the sequences $\{r_0+r_1+\hdots + r_{n-1}\}$
%and $\{\la (n-1)+\nu\}$ for big $n$,
%remembering that $r_k$ is non-decreasing,
%we now obtain the following
%\begin{prop}\label{prop:Vnumber}
%$r_k\leq \la$ for all $k$ and
%there exists a number $N$ such that $r_{N+n}=r_N$ for all $n\geq 0$.
%\end{prop}
%
We get the following result as a direct consequence
of the results of this part.
%of the Proposition above, Theorem~\ref{thm:main2},
%Corollary~\ref{cor:v=v} and Theorem~\ref{thm:fullKM}.
\begin{prop}
Let $p$ be a semi-regular prime.
Then all three of the Kervaire-Murthy
conjectures hold if and only if $\la=r(p)$. Moreover,
if $\la=r(p)$
then $\nu=r(p)$.
\end{prop}

{\bf The case $\la=r_1$}

In \cite{U2} Ullom uses the following somewhat technical assumption on $p$.
\begin{ass}\label{ass:Ullomcond}
The Iwasawa
invariants $\la_{1-i}$ satisfy\label{la_1-i}
$
1\leq\la_{1-i}\leq p-1
$
\end{ass}
We refer you to \cite{I} for notation.
S. Ullom proves that if Assumption~\ref{ass:Ullomcond} holds then, for even $i$,
\begin{equation}\label{eq:Ullom}
e_i V_n\cong \frac{\Z}{p^n \Z}\oplus (\frac{\Z}{p^{n-1} \Z})^{\la_{1-i}-1}.
\end{equation}
This yields, under the same assumption, that
\begin{equation}\label{eq:Ullom2}
V_n^+\cong (\frac{\Z}{p^n \Z})^{r(p)} \oplus
(\frac{\Z}{p^{n-1}\Z})^{\la-r(p)},
\end{equation}
where
\[
\la=\sum_{i=1,\, i \text{ even}}^{r(p)} \la_{1-i}
\]
Hence,
%under the assumption $\la_{1-i}=0$ or $1$
when $\la=r(p)$
we get
%\ref{eq:conj2}.
that $V_n=(\Z/p^n \Z)^ {r(p)}$ as predicted by Kervaire and Murthy.
Note however, that if $\la>r(p)$, then Kervaire and Murthy's conjecture fails.
%conjecture~\ref{eq:conj2} is false.
We will discuss some consequenses of Assumption 2.
Since $V_n^+$ is a quotient of
$\V_n^+$ applying this to $n=n_0+1$ yields
\[
r_0 + \la n_0 \leq r_0+r_1+\hdots+r_{n_0} \leq r_0+n_0 r_{n_0} \leq r_0+n_0\la
.
\]
This obviously implies that $r_k=\la$ for all $k=1,2,\hdots$ because
$\{r_i\}$ is a non-decreasing sequence
bounded by $\la$ by Proposition 2.11 and hence we get the following
\begin{lemma}\label{lemma:la=rn}
When Assumtion~\ref{ass:Ullomcond}
holds $r_k=\la$ for all $k=1,2,\hdots$.
\end{lemma}
The following theorem is now immediate.
\begin{thm}
If Assumtion~\ref{ass:Ullomcond} holds, then $\V_n^+=V_n^+$.
\end{thm}
Keeping Theorem~\ref{thm:fullKM} in mind we immediately get the following corollary.
\begin{cor}
When Assumtion~\ref{ass:Ullomcond} holds,
$$\cl^{(p)} \Q(\z_{n-1})\cong (\Z/p^{n}\Z)^{r(p)} \oplus (\Z/p^{n-1}\Z)^{\la-r(p)}$$
and
$\la=r_1$ and $\nu=r_0=r(p)$.
\end{cor}

\subsection{An application to units in $\Z[\z_n]$} \mbox{}
%=\subsection{Missed Places}
%
%
%
%
%

The techniques we have developed also lead to some conclusions about the group of units
in $\Z[\z_n]^*$. From the previous results we know that
\[
\V_{n+1}^+=\frac{\tilde{D}_{n+1}^{*+}}{g_{n+1}(U_{n,1})}\cong \
\frac{\tilde{D}_{n+1}^{*+}}{\frac{U_{n,2}}{U_{n,p^{n+1}-1}}}
\]
Let  $s_{n,p^{n+1}-1}=|U_{n,1}/U_{n-1,p^{n+1}-1}|$. A naive first guess would be that
$s_{n,p^{n+1}-1}=\frac{p^{n+1}-1-2}{2}=\frac{p^{n+1}-3}{2}$ which is the maximal value of this number.
Incidentally, this maximal value equals $|\tilde{D}_{n+1}^{*+}|$. In this case we say that
$U_{n,1}/U_{n,p^{n+1}-1}$ is full, but this happens if and only if $p$ is a regular prime.
In other words $\V_{n+1}^+$ is trivial if and only if $p$ is regular.
This fact is by the way proved directly in \cite{H}.
For non-regular (but as before semi-regular) primes what happens is that there are
``missed places'' in $U_{n,1}/U_{n,p^{n+1}-1}$. We define $2k$ as a missed place (at level $n$) if
$U_{n,2k}/U_{n,2k+2}$ is trivial.
Lemma~\ref{lemma:p^k-1=p^k+1} reads
$U_{n,p^{n+1}-1}=U_{n,p^{n+1}+1}$ and hence
provides an instant example of a missed place, namely $p^{n+1}-1$. It follows from our theory
that every missed place corresponds to a non-trivial element of $\V_{n+1}^+$.
Recall that $\Z[\z_{n-1}]^*$ is identified with its image in $A_n$. We will now prove that the map
$g_n:\Z[\z_{n-1}]^*\to D_n^*$ respects the filtrations $\la_{n-1}^k$ and $(x-1)^k$.
\begin{prop}
Let $1\leq s \leq p^n-1$ and $\e\in \Z[\z_{n-1}]^*$.
Then $\e\in U_{n-1,s}$ if and only if $g_n(\e)\in D_{n,(s)}^*$.
\end{prop}
Using this Proposition we see that an element of $D_{n+1,(2s)}^{*+}$ which is non-trivial in $\V_{n+1}^+$
corresponds to a missed place $2s$ at level $n$.
\begin{proof}
To show that
$g_n(\e)\in D_{n,(s)}^*$ implies $\e\in U_{n-1,s}$ we can use the same technique as in the proof of
Theorem I.2.7 in \cite{ST98} (also see Lemma~\ref{lemma:gker}).
For the other direction, first
note that is $s\leq p^n-p^{n-1}$ the statement follows directly from the commutativity of the diagram
\begin{equation}\label{eq:Z+Amodpcomm}
\xymatrix@=40pt{
A_n^* \ar[r] \ar[d]^{\m p}  & \Z[\z_{n-1}]^{*+}  \ar[d]^{\m p}\\
  D_{n}^*\cong \big(\frac{\F_p[x]}{(x-1)^{p^n-1}}\big)^* \ar[r] & \big(\frac{\F_p[x]}{(x-1)^{p^n-p^{n-1}}}\big)^* }
\end{equation}
What is left to prove is that $\e\in U_{n-1,s}$ implies $g_n(\e)\in D_{n,(s)}^*$ also for
$p^n-p^{n-1}\leq s \leq p^n-1$. For technical reason we will prove that
if $\e\in U_{n-1,p^n-p^k+r}$ for some $1\leq k \leq n-1$ and
$0\leq r \leq p^k-p^{k-1}$ then $g_n(\e)\in D_{n,(p^n-p^k+r)}^*$.
Note that $\e\in U_{n-1,p^n-1}$ is equivalent to $g_n(\e)=1\in D_{n}^*$ by Lemma~\ref{lemma:gker}.
Suppose $\e=1+t\la_{n-1}^{p^n-p^k+r}$ for some $t\in \Z[\z_{n-1}]$.
By
%Lemma 4.11 of \cite{preprint-2-01}
Lemma~\ref{lemma:normorder}
we get  $N_{n-1}(\e)=1+t'p(x-1)^{p^{n-1}-p^{k-1}}$ for some $t' \in A_{n-1}$. In $A_n$,
\[
p=\frac{x^{p^{n}}-1}{x^{p^{n-1}}-1}+t(x)\frac{x^{p^{n-1}}-1}{x-1}
\]
for some polynomial $t(x)$. In $A_n$ consider the element
%$p(t(x-1)^{p^{n-1}-p^k+r}-t'(x-1)^{p^{n-1}-p^{k-1}})$ in $A_n$.
%By using the expression for $p$ above, we get
\begin{equation*}
\begin{split}
p(t(x-1)^{p^{n-1}-p^k+r}-t'(x-1)^{p^{n-1}-p^{k-1}})=\\
=\big(\frac{x^{p^{n}}-1}{x^{p^{n-1}}-1}+t(x)\frac{x^{p^{n-1}}-1}{x-1}\big)
(t(x-1)^{p^{n-1}-p^k+r}-t'(x-1)^{p^{n-1}-p^{k-1}}).
\end{split}
\end{equation*}
By computing the right hand side and re-arrange the terms we get
\begin{equation*}
\begin{split}
f:=tp(x-1)^{p^{n-1}-p^k+r}-\big(t(x-1)^{p^{n-1}-p^k+r}-t'(x-1)^{p^{n-1}-p^{k-1}}\big)\frac{x^{p^{n}}-1}{x^{p^{n-1}}-1}=\\
=t'(x-1)^{p^{n-1}-p^{k-1}}- b(x)\frac{x^{p^{n-1}}-1}{x-1}.
\end{split}
\end{equation*}
Using the two representations of $f$ we see that $i_n(1+f)=\e$ and $j_n(1+f)=N_{n-1}(\e)$ so
$1+f$ represents $(\e,N_{n-1}(\e))$ (which represents $\e$ under our usual identification) in $A_n$.
Since $\leq p^k-p^{k-1}$ we now get
$g_n(1+f)\equiv 1 \m (x-1)^{p^{n-1}-p^k+r}$ in $D_n$ as asserted.
\end{proof}

%By Lemma 5.2 of \cite{preprint-2-01} we have that for all $1\leq 2s\leq p^{n+1}-1$,
%$\e\in U_{n,2s}$ if and only if an only if
%$g_{n+1}(\e)\in D_{n+1}^{*+}(2s)$, where  $D_{n+1}^{*}(k)=\{ a:a=1\m(x_{n+1}-1)^k\}$.
Theorem~\ref{thm:Vn-semireg} and its proof
now give us specific information about the missed places which we will formulate in a Theorem below.
We start with a simple lemma.
\begin{lemma}
Let $1\leq s \leq n+1$ and $1\leq k < s$. Then
$p^s-p^k$ is a missed place at level $n$ if and only if $s=n+1$ and $k=1$.
\end{lemma}
\begin{proof}
Let $\eta:=\z_n^{(p^{n+1}+1)/2}$. Then $\eta^2=\z_n$ and $c(\eta)=\eta^{-1}$. Define
\[
\e:=\frac{\eta^{p^s+p^k}-\eta^{-(p^s+p^k)}}{\eta^{p^k}-\eta^{-(p^k)}}.
\]
Clearly, $\e$ is real and since
\[
\e=\eta^{-p^s}\frac{\z_n^{p^s+p^k}-1}{\z_n^{p^k}-1},
\]
$\e$ is a unit.
By a calculation one can show that $\e\in U_{n,p^s-p^k}\setminus U_{n,p^s-p^k+2}$.
\end{proof}
Define for $k=0,1,\hdots$ the $k$-strip as the numbers $p^k+1,p^k+3,\hdots,p^{k+1}-1$.
\begin{thm}
At level $n$ we have the following
\begin{enumerate}
\item Let $0 \leq k \leq n$. In the $k$-strip there are exactly $r_k$ missed places.
\item The missed places in the 0-strip are in one to one correspondence with the numbers
$2i_1,\hdots,2i_{r_0}$ such that the numerator of the Bernoulli-number $B_{2i_k}$ (in reduced form)
is divisible by $p$.
\item Suppose $i_1,\hdots,i_{r_k}$ are the missed places in the $k$-strip. Then $pi_1,\hdots,pi_{r_k}$
are missed places in the $k+1$ strip. The other $r_{k+1}-r_{k}$ missed places in the $k+1$ strip
are not divisible by $p$.
\end{enumerate}
\end{thm}
\begin{proof}
We know from  Proposition~\ref{prop:V1-generators}
%Proposition 4.6 of \cite{preprint-2-01}
that we have $r_0$ missed  places
in the 0-strip at level 0 and that they correspond exactly to the indexes of the relevant Bernoulli
numbers.
As in Proposition~\ref{prop:Dn-Zn-generators}
%Proposition 4.7 of \cite{preprint-2-01}
an induction argument
using the map $\pi_n$ to lift the generators of $\V_{n-1}^+$ to $\V_{n}^+$
show that we have $r_0$ missed places in the 0-strip at every level and that
a missed place $k$ at level $n-1$ lift to missed places $k$ and $pk$ at level $n$.
%This proves {\it 2}.
%To prove {\it 1} we now only need prove that the
What is left to prove is that the
``new'' missed places we get when we go
from level $n-1$ to $n$ all end up in the $n$-strip
and that no ``new'' missed places are divisible by $p$.
First, $p^n-1$ can not be a
missed place (at level $n$) by the lemma above.
%since the unit $\e_n$ from the proof of Proposition~\ref{Vn=Vn}
%satisfies $\e_n\in U_{n,p^n-1}\setminus U_{n,p^n+1}$.
It follows from our theory that the ``new'' missed
places correspond to the generators of $\V_{n+1}^+$ of exponent $p$.
We need to show that each such generators
$a_l$, $l=1,\ldots,r_{n-1}-r_{n-2}$, belong to $D_{n+1,(p^n+1)}^{*+}$. Suppose for a contradiction
that $a_l=1+t(x_{n+1}-1)^s$, $t\neq 0$, $s< p^n-1$, is a ``new'' generator. Then
$\pi_{n+1}(a_l)=1+t(x_{n}-1)^s$ is neccesarily trivial in $\V_{n}^+$ but not in $D_n^{*+}$.
%But the map $\bar{j_{n+1}}$ which induces $\pi_{n+1}$ maps $1+t(x_{n+1}-1)^s$ to $1+t(x_{n}-1)^s$ which is non-trivial in $\tilde{D}_{n}^{*+}$.
Hence $\pi_{n+1}(a_l)=g_{n}(\e)$ for some $\e\in \Z[\z_{n-1}]^*$. Since the usual norm map
$\tilde{N}_{n,1}$ is surjective (when $p$ is semi-regular) and by commutativity of
diagram~\ref{eq:ZDcomm}
%diagram 4.1 of  \cite{preprint-2-01}
we then get
%$1+t(x_{n+1}-1)^s=g_{n}(\e')$
$a_lg_{n+1}(\e')^{-1}=b$ for some
$\e'\in \Z[\z_{n}]^*$
and
$b\in \ker\{\tilde{D}_{n+1}^{*+}\rightarrow \tilde{D}_{n}^{*+}\}=\tilde{D}_{n+1}^{*+}(p^n-1)$.
Since $p^n-1$ is not a missed place,
$b=g_{n+1}(\e'')$ for some
some $\e''\in \Z[\z_{n}]^*$. But this
means $a_l$ is trivial in $\V_{n+1}^+$
which is a contradiction. We conclude that $a_l\in D_{n+1,(p^n+1)}^{*+}$.

To prove no ``new'' missed places are divisible by $p$
%we use the map $\alpha_n$ to see that a missed place $k$ at level $n-1$ lifts to a
%missed place $pk$ at level $n$. To prove the rest of {\it 3} it is enough to prove that no ``new''
%missed places are divisible by $p$ (since the rest follows inductively).
we need
to show that if $a_l\in D_{n+1,(s)}^{*+}\setminus D_{n+1,(s+2)}^{*+}$
is a ``new'' generator of $\V_{n+1}^+$, then $p$
does not divide $s$.
Now, a generator can always be chosen of the form $1+(x_{n+1}-1)^{s}$.
Then an element of the form $1+(x_{n+1}-1)^{pk}$, with
$k\not\in \{i_1,\hdots,i_{r_{n-1}}\}$ cannot be a missed place.
This follows from the fact that if $k$ is not a missed place, then $1+(x_{n}-1)^{k}$ is trivial in $\V_n^+$
and since $\alpha_n$ is injective,
$1+(x_{n+1}-1)^{pk}=\alpha_n(1+(x_{n}-1)^{k})$ is also trivial in $\V_{n+1}^+$.
%Suppose the new generator $a_l\in D_{n+1}^{*+}(pt)\setminus D_{n+1}^{*+}(pt+2)$
%is arbitrary. Then
%$a_l=(1+a(x_{n+1}-1)^{pt})b$ for some $b\in D_{n+1}^{*+}(pt+2)$. But since $a_l$ is a new generator,
%$a_l^p$ is trivial and hence $b^p$ is trivial. $b$ can not be one of the new generators and hence
%$b$ is the $p$-th power of an old generator.
%%This means
%%\[
%%a_l=1+a_{pt}(x_{n+1}-1)^{pt}+a_{p(t+1)}(x_{n+1}-1)^{p(t+1)}+\hdots=
%%\alpha_n(1+a_{pt}(x_{n}-1)^{t}+a_{p(t+1)}(x_{n}-1)^{(t+1)}+\hdots).
%%\]
%As in the case $(1+a(x_{n+1}-1)^{pt})$ this means $a_l$ is trivial which is a contradiction.
\end{proof}

\end{document}